\documentclass{amsart}
\usepackage{amssymb}
\usepackage{xcolor,graphicx}

\theoremstyle{plain}
\newtheorem{theorem}{Theorem}
\newtheorem{lemma}[theorem]{Lemma}

\theoremstyle{definition}

\theoremstyle{remark}
\newtheorem{remark}[theorem]{Remark}
\newtheorem{example}[theorem]{Example}

\providecommand{\abs}[1]{\lvert#1\rvert}
\providecommand{\Abs}[1]{\Bigl\lvert#1\Bigr\rvert}
\providecommand{\norm}[1]{\lVert#1\rVert}

\newcommand\independent{\protect\mathpalette{\protect\independenT}{\perp}}
\def\independenT#1#2{\mathrel{\rlap{$#1#2$}\mkern2mu{#1#2}}}

\begin{document}

\title[Knockoffs]{Knockoffs for exchangeable\\categorical covariates}

\author{Emanuela Dreassi}
\address{Emanuela Dreassi, Dipartimento di Statistica, Informatica, Applicazioni ``G. Parenti'', Universit\`a di Firenze, viale Morgagni 59, 50134 Firenze, Italy}
\email{emanuela.dreassi@unifi.it}
\author{Luca Pratelli}
\address{Luca Pratelli, Accademia Navale, viale Italia 72, 57100 Livorno,
Italy} \email{luca{\_}pratelli@marina.difesa.it}
\author{Pietro Rigo}
\address{Pietro Rigo (corresponding author), Dipartimento di Scienze Statistiche ``P. Fortunati'', Universit\`a di Bologna, via delle Belle Arti 41, 40126 Bologna, Italy}
\email{pietro.rigo@unibo.it}

\begin{abstract}
Let $X=(X_1,\ldots,X_p)$ be a $p$-variate random vector and $F$ a finite set. In a number of applications, mainly in genetics, it turns out that $X_i\in F$ for each $i=1,\ldots,p$. Despite the latter fact, to obtain a knockoff $\widetilde{X}$ (in the sense of \cite{CFJL18}), $X$ is usually modeled as an absolutely continuous random vector. While comprehensible from the point of view of applications, this approximate procedure does not make sense theoretically, since $X$ is supported by the finite set $F^p$. In this paper, explicit formulae for the joint distribution of $(X,\widetilde{X})$ are provided when $P(X\in F^p)=1$ and $X$ is exchangeable or partially exchangeable. In fact, when $X_i\in F$ for all $i$, there seem to be various reasons for assuming $X$ exchangeable or partially exchangeable. The robustness of $\widetilde{X}$, with respect to the de Finetti's measure $\pi$ of $X$, is investigated as well. Let $\mathcal{L}_\pi(\widetilde{X}\mid X=x)$ denote the conditional distribution of $\widetilde{X}$, given $X=x$, when the de Finetti's measure is $\pi$. It is shown that
$$\norm{\mathcal{L}_{\pi_1}(\widetilde{X}\mid X=x)-\mathcal{L}_{\pi_2}(\widetilde{X}\mid X=x)}\le c(x)\,\norm{\pi_1-\pi_2}$$
where $\norm{\cdot}$ is total variation distance and $c(x)$ a suitable constant. Finally, a numerical experiment is performed. Overall, the knockoffs of this paper outperform the alternatives (i.e., the knockoffs obtained by giving $X$ an absolutely continuous distribution) as regards the false discovery rate but are slightly weaker in terms of power.
\end{abstract}

\maketitle

\section{Introduction}\label{intro}

\noindent The knockoff procedure (KP) is a powerful and effective tool for making variable selection in regression problems. It has been introduced by Barber and Candes in \cite{BC15} and subsequently developed by various authors. Without any claim of being exhaustive, a list of references is: \cite{BCS20, BCJW21, JSPI, CFJL18, EJS, GGZ, KATSAB, RC, RWC, RSC, SSC19, WSBBC}.

\medskip

\noindent To formalize, let $Y$ be the response variable and
$$X=(X_1,\ldots,X_p)$$
the vector containing the $p$ covariates $X_1,\ldots,X_p$. Moreover, for any random variables $U$ and $V$, let $\mathcal{L}(U)$ denote the probability distribution of $U$ and $\mathcal{L}(U\mid V)$ the conditional distribution of $U$ given $V$. The naive idea of KP is to build an auxiliary vector
$$\widetilde{X}=(\widetilde{X}_1,\ldots,\widetilde{X}_p),$$
called a {\em knockoff copy of} $X$, which is able to capture the connections among $X_1,\ldots,X_p$. Once $\widetilde{X}$ is given, each $X_i$ is selected/discarded based on the comparison between it and $\widetilde{X}_i$. Roughly speaking, $\widetilde{X}_i$ is a control for $X_i$, and $X_i$ is selected if it appears to be considerably more associated with $Y$ than its knockoff copy $\widetilde{X}_i$. How to build $\widetilde{X}$ is discussed in Section \ref{vh87}. Here, we just recall three main features of KP:

\begin{itemize}

\item KP allows to take the false discovery rate under control;

\item For KP to apply, the statistician must assign $\mathcal{L}(X)$ but is not forced to specify $\mathcal{L}(Y\mid X)$;

\item For KP to have a good performance in terms of power, it is desirable that $X$ and $\widetilde{X}$ are ``as independent as possible".

\end{itemize}

\medskip

\subsection{Motivations}
In a number of applications, there is a {\em finite} set $F$ such that
\begin{gather}\label{i9p}
X_i\in F\quad\quad\text{for each }i=1,\ldots,p.
\end{gather}
In genetics, for instance, it is quite usual that $F=\{0,1\}$ and $X_i=1$ or $X_i=0$ according to whether the $i$-th gene is present or not. Or else, $F=\{0,1,2\}$ where $X_i=0$ and $X_i=1$ have the previous interpretation while $X_i=2$ means that the absence/presence of the $i$-th gene cannot be established. Condition \eqref{i9p} is also quite recurrent in causal inference, mainly in connection with enviromental or biomedical applications; see  e.g. \cite{CAUS1, CAUS2} and references therein.

\medskip

\noindent To our knowledge, from the theoretical point of view, condition \eqref{i9p} has been almost neglected to date (the only exceptions known to us are \cite{DAIBARB, EJS, STATMED}). More precisely, there is a plenty of papers applying KP to categorical covariates. Despite of \eqref{i9p}, however, $X$ is given an absolutely continuous distribution. This choice of $\mathcal{L}(X)$ is justified by the need of the application and/or by the difficulty of constructing a reasonable knockoff. It should be added that, in most applications, KP exhibits a good behavior even if $X$ is given an absolutely continuous distribution. Nevertheless, under condition \eqref{i9p}, taking $\mathcal{L}(X)$ absolutely continuous does not make sense. It is merely an approximate solution which strongly contrast with a basic feature of the data. Note also that, since $\mathcal{L}(X)=\mathcal{L}(\widetilde{X})$, condition \eqref{i9p} yields
$$\widetilde{X}_i\in F\quad\quad\text{a.s. for each }i=1,\ldots,p.$$

\medskip

\noindent This paper investigates KP, under condition \eqref{i9p}, from the theoretical point of view. Our main goal is to understand what happens when condition \eqref{i9p} is taken seriously into account and $X$ is given a distribution such that
$$P(X\in F^p)=1.$$
We focus on the special case where $X$ is (infinitely) exchangeable or (infinitely) partially exchangeable. These assumptions, like any other, have merits and drawbacks. Under condition \eqref{i9p}, however, they seem to be reasonable.

\medskip

\noindent To support the latter claim, we list some consequences of exchangeability and partial exchangeability. Before doing this, we briefly introduce some notation. Without loss of generality, we let
$$F=\bigl\{0,1,\ldots,m\bigr\}\quad\quad\text{for some integer }m\ge 1.$$
If $x$ is a point of $F^p$ or $\mathbb{R}^p$, we denote by $x_i$ the $i$-th coordinate of $x$. Moreover, for all $j\in F$ and all $x\in F^p$, we write
$$n_j=n_j(x)=\sum_{i=1}^p\textbf{1}(x_i=j)$$
where $\textbf{1}(x_i=j)$ is the indicator of the event $\{x_i=j\}$.

\medskip

\noindent Exchangeability of $X$ can be motivated as follows.

\medskip

\begin{itemize}

\item[(i)] $P(X_i=j)$ does not depend on $i$. This makes sense when the statistician is very ignorant about the effect of the covariates, and he/she does not want favour one over the others.

\medskip

\item[(ii)] $P(X=x)$ depends only on $n_0(x),\ldots,n_m(x)$. For  large $p$, this is fundamental from the numerical point of view. Otherwise, $P(X=x)$ is to be storaged for each $x\in F^p$.

\medskip

\item[(iii)] The order in which $X_1,\ldots,X_p$ are observed provides no information, and this is sound in several applications.

\medskip

\item[(iv)] To build a knockoff $\widetilde{X}$ is straightforward and $\mathcal{L}(\widetilde{X}\mid X)$ can be written in closed form; see below.

\end{itemize}

\medskip

\noindent As to partial exchangeability, suppose that $X$ is partitioned into subvectors, say $X'=(X_1,\ldots,X_k)$ and $X''=(X_{k+1},\ldots,X_p)$, which are to be modeled differently; see e.g. \cite{KATSAB}. For instance, the elements of $X'$ are more reliable than those of $X''$. Or else, there is a latent factor which affects $X'$ and $X''$ in a different way. At the same time, the statistician is quite ignorant about the covariates within $X'$ and is not able to discriminate among them. The same happens for $X''$. In situations of this type, it may be reasonable to assume $X_1,\ldots,X_p$ conditionally independent given a random vector
$$Z=(U_1,\ldots,U_m,V_1,\ldots,V_m).$$
Moreover, letting $U_0=1-\sum_{j=1}^mU_j$ and $V_0=1-\sum_{j=1}^mV_j$, one also assumes
$$P(X_i=j\mid Z)=U_j\quad\text{or}\quad P(X_i=j\mid Z)=V_j$$
according to whether $i\le k$ or $i>k$. This means that $X$ is partially exchangeable. Note that, since $X'$ and $X''$ are marginally exchangeable, points (i)-(ii)-(iii) still apply with reference to each single subvector. In turn, point (iv) holds by the results of \cite{EJS}. See Sections \ref{prel} and \ref{bv5r4}.

\medskip

\noindent In a nutshell, this paper investigates KP when condition \eqref{i9p} holds and $X$ is given an exchangeable or partially exchangeable distribution. Indeed, under \eqref{i9p}, this modeling of $X$ can be supported by the previous arguments.

\medskip

\noindent In general, unless $X$ satisfies suitable assumptions, to build a knockoff is very hard. Hence, Remark (iv) is fundamental. We spend some more on it when $F=\bigl\{0,1\bigr\}$ and $X$ is (infinitely) exchangeable. In this case, de Finetti's theorem yields
$$P(X=x)=\int_{[0,1]}u^{n_1}(1-u)^{n_0}\,\pi(du),\quad\quad x\in\{0,1\}^p,$$
for some (unique) probability measure $\pi$ on $[0,1]$. Such a $\pi$ is usually called {\em the de Finetti's measure of} $X$. Even if the Bayesian interpretation is not mandatory, in the sequel, $\pi$ is also referred to as {\em ``the prior"}.

\medskip

\noindent If the statistician agrees to assume $X$ exchangeable, the only remaining task is to select the prior $\pi$. After doing this, to obtain a knockoff, it suffices to let
\begin{gather}\label{p9ij5x}
P\bigl(X=x,\,\widetilde{X}=\widetilde{x}\bigr)=\int_{[0,1]} u^{n_1+\widetilde{n}_1}(1-u)^{n_0+\widetilde{n}_0}\,\pi(du)
\end{gather}
for all $x,\,\widetilde{x}\in\{0,1\}^p$, where $\widetilde{n}_j=n_j(\widetilde{x})=\sum_{i=1}^p\textbf{1}(\widetilde{x}_i=j)$. In fact, by Theorem \ref{gy78n3} below, any $p$-variate random vector $\widetilde{X}$ satisfying equation \eqref{p9ij5x} is a knockoff.

\medskip

\noindent Note also that equation \eqref{p9ij5x} yields
$$P\bigl(\widetilde{X}=\widetilde{x}\mid X=x\bigr)=\frac{P\bigl(X=x,\,\widetilde{X}=\widetilde{x}\bigr)}{P(X=x)}=\frac{\int_{[0,1]} u^{n_1+\widetilde{n}_1}(1-u)^{n_0+\widetilde{n}_0}\,\pi(du)}{\int_{[0,1]} u^{n_1}(1-u)^{n_0}\,\pi(du)}$$
provided $P(X=x)>0$. Thus, $\mathcal{L}(\widetilde{X}\mid X)$ can be written in closed form. In applications, this is very useful. In fact, after observing $X=x$, the values $\widetilde{x}$ of $\widetilde{X}$ (required to implement KP) can be sampled directly from $\mathcal{L}(\widetilde{X}\mid X=x)$.

\medskip

\subsection{Our contribution} The layout of the paper is as follows. Section \ref{prel} provides some preliminaries. Sections \ref{b76t}-\ref{rob1} are the core of the paper and include the main results. Explicit formulae for $\mathcal{L}(X,\widetilde{X})$ are given when $P(X\in F^p)=1$ and $X$ is exchangeable or partially exchangeable. For the sake of simplicity, in the partially exchangeable case, we focus on $F=\{0,1\}$ only, but our argument applies to $F=\bigl\{0,1,\ldots,m\bigr\}$ for any $m\ge 1$. To evaluate $\mathcal{L}(X,\widetilde{X})$, since $X$ is exchangeable or partially exchangeable, the de Finetti's measure $\pi$ of $X$ is to be selected. Various choices of $\pi$, both discrete and diffuse, are taken into account. Discrete $\pi$'s have an advantage: they admit a natural interpretation which makes their choice easier. To illustrate, suppose $X$ is exchangeable, $F=\{0,1\}$, and one decides to use a discrete $\pi$. Then, to specify $\pi$, it suffices to answer the question: ``what about $P(n_1=k)$, i.e., what is the probability of exactly $k$ successes in $p$ trials ?" Answering this question, for each $k\in\{0,1,\ldots,p\}$, amounts to selecting $\pi$. Some other results, contained in Section \ref{rob1}, concern the robustness of $\widetilde{X}$ with respect to the choice of $\pi$. Let $\mathcal{L}_\pi(\widetilde{X}\mid X=x)$ denote the conditional distribution of $\widetilde{X}$, given $X=x$, when the de Finetti's measure is $\pi$. Then, among other things, it is shown that
$$\norm{\mathcal{L}_{\pi_1}(\widetilde{X}\mid X=x)-\mathcal{L}_{\pi_2}(\widetilde{X}\mid X=x)}\le c(x)\,\norm{\pi_1-\pi_2}$$
where $\norm{\cdot}$ is total variation distance and $c(x)$ a suitable constant. Finally, Section \ref{b5f9k} includes some numerical illustrations for simulated and real data. Overall, the knockoffs of this paper outperform the alternatives (i.e., the knockoffs obtained by giving $X$ an absolutely continuous distribution) as regards the false discovery rate, but not necessarily in terms of power.

\medskip

\section{Preliminaries}\label{prel}

\medskip

\subsection{Exchangeability}\label{nd56c} An (infinite) sequence $(Z_n:n\ge 1)$ of real random variables is {\em exchangeable} if
$$(Z_{\sigma_1},\ldots,Z_{\sigma_n})\sim (Z_1,\ldots,Z_n)$$
for all $n\ge 1$ and all permutations $(\sigma_1,\ldots,\sigma_n)$ of $(1,\ldots,n)$. An i.i.d. sequence is obviously exchangeable, while the converse is not true (just take $Z_n=V$ for each $n$, where $V$ is any non degenerate random variable). However, by de Finetti's theorem, $(Z_n)$ is exchangeable if and only if there is a random probability measure $\mu$ on $\mathbb{R}$ such that
$$(Z_n)\text{ is conditionally i.i.d., given }\mu,\text{ and }P(Z_1\in\cdot\mid\mu)=\mu(\cdot)\text{ a.s.}$$
Equivalently,
\begin{gather}\label{pres}
P\bigl(Z_1\in A_1,\ldots,Z_n\in A_n)=E\bigl\{\mu(A_1)\ldots\mu(A_n)\bigr\}
\end{gather}
for all $n\ge 1$ and all $A_1,\ldots,A_n\in\mathcal{B}(\mathbb{R})$.

\medskip

\noindent The probability distribution $\pi$ of $\mu$ is said to be the {\em de Finetti's measure} of the sequence $(Z_n)$. In general, $\pi$ is a probability on $\mathcal{P}(\mathbb{R})$, where $\mathcal{P}(\mathbb{R})$ is the set of all probability measures on $\mathbb{R}$. If the $Z_n$ take values in a finite set $F$, however, $\pi$ can be regarded as a probability on a suitable finite dimensional space. For instance, if $F=\{0,1\}$, $\pi$ can be viewed as a probability on $[0,1]$. More generally, if $F=\bigl\{0,1,\ldots,m\bigr\}$, $\pi$ can be seen as a probability on the simplex
$$S=\bigl\{u\in\mathbb{R}^m:u_j\ge 0\text{ for all }j\text{ and }\sum_{j=1}^mu_j\le 1\bigr\}.$$

\medskip

\noindent Finally, a $n$-variate random variable $Z=(Z_1,\ldots,Z_n)$ is exchangeable provided $(Z_{\sigma_1},\ldots,Z_{\sigma_n})\sim (Z_1,\ldots,Z_n)$ for all permutations $(\sigma_1,\ldots,\sigma_n)$ of $(1,\ldots,n)$. If condition \eqref{pres} holds then $Z$ is exchangeable, but the converse is not true. To avoid misunderstandings, in this paper, a $n$-variate random variable $Z$ is said to be {\em infinitely exchangeable} if condition \eqref{pres} holds for some random probability measure $\mu$.

\medskip

\subsection{Partial exchangeability} Let $\{C_1,\ldots,C_m\}$ be a partition of $\{1,2,\ldots\}$. Define a function $f:\{1,2,\ldots\}\rightarrow\{1,\ldots,m\}$ as $f(i)=j$ if $i\in C_j$. Then, $(Z_n)$ is {\em partially exchangeable} if $(Z_{\sigma_1},\ldots,Z_{\sigma_n})\sim (Z_1,\ldots,Z_n)$ for all $n\ge 1$ and all permutations $(\sigma_1,\ldots,\sigma_n)$ of $(1,\ldots,n)$ such that
$$f(\sigma_i)=f(i)\quad\quad\text{for all }i=1,\ldots,n.$$
Thus, the probability distribution of $(Z_1,\ldots,Z_n)$ is not necessarily invariant under all permutations, but only under those which swap elements in the same equivalence class. This weaker assumption allows to model much more situations than exchangeability. Moreover, a version of de Finetti's theorem is still available; see \cite{ALDOUS}. Suppose $C_j$ is infinite for each $j$. Then, $(Z_n)$ is partially exchangeable if and only if there are $m$ random probability measures $\mu_1,\ldots,\mu_m$ on $\mathbb{R}$ such that
\begin{gather}\label{7yhw}
P\bigl(Z_1\in A_1,\ldots,Z_n\in A_n)=E\left\{\prod_{i=1}^n\mu_{f(i)}(A_i)\right\}
\end{gather}
for all $n\ge 1$ and all $A_1,\ldots,A_n\in\mathcal{B}(\mathbb{R})$. Finally, let $Z=(Z_1,\ldots,Z_n)$ be a $n$-variate random variable. Suppose that $\{1,\ldots,n\}$ is partitioned into the sets $C_1,\ldots,C_m$. Then, we say that $Z$ is {\em infinitely partially exchangeable} if condition \eqref{7yhw} holds for some random probability measures $\mu_1,\ldots,\mu_m$.

\medskip

\subsection{Knockoffs and conditional independence}\label{vh87} For $i=1,\ldots,p$, define a function $f_i:\mathbb{R}^{2p}\rightarrow\mathbb{R}^{2p}$ as
$$f_i(x)=(x_1,\ldots,x_{p+i},\ldots,x_p,x_{p+1}\ldots,x_i,\ldots,x_{2p}),\quad x\in\mathbb{R}^{2p}.$$
Thus, $f_i(x)$ is the point of $\mathbb{R}^{2p}$ obtained by swapping $x_i$ with $x_{p+i}$ and leaving all other coordinates of $x$ fixed. A $p$-variate random vector $\widetilde{X}$ is a knockoff copy of $X$, or merely a knockoff, if
$$f_i(X,\widetilde{X})\sim (X,\widetilde{X})\text{ for }i=1,\ldots,p\quad\text{and}\quad\widetilde{X}\independent Y\mid X.$$
Note that a knockoff always exists. It suffices to let $\widetilde{X}=X$. This trivial knockoff, however, is unsuitable for KP to apply.

\medskip

\noindent As noted in Section \ref{intro}, to implement KP, one must select $\mathcal{L}(X)$. However, if $\mathcal{L}(X)$ is arbitrary, the construction of a (reasonable) knockoff $\widetilde{X}$ is usually hard. There are some universal algorithms, such as the {\em Sequential
Conditional Independent Pairs} \cite{CFJL18} and the {\em Metropolized Knockoff Sampler} \cite{BCJW21},
which are virtually able to cover any choice of $\mathcal{L}(X)$. Nevertheless, these algorithms are quite difficult to use and
do not provide a closed formula for $\widetilde{X}$. More importantly, they are computationally
intractable when $\mathcal{L}(X)$ is complex; see \cite{BCJW21} and \cite{GGZ}.

\medskip

\noindent Leaving algorithms aside, to our knowledge, the most general condition for obtaining a closed formula for $\widetilde{X}$ is
\begin{gather}\label{m9i}
X_1,\ldots,X_p\,\text{ conditionally independent given some random element }Z;
\end{gather}
see \cite{JSPI} and \cite{EJS}. We next briefly describe how to build $\widetilde{X}$ under \eqref{m9i}.

\medskip

\noindent Let $\Theta$ denote the set where $Z$ takes values, $\pi$ the probability distribution of $Z$, and
$$P_i(A\mid\theta)=P(X_i\in A\mid Z=\theta)\quad\quad\text{for }i=1,\ldots,p,\,\theta\in\Theta\text{ and }A\in\mathcal{B}(\mathbb{R}).$$
Because of \eqref{m9i},
$$P(X_1\in A_1,\ldots,X_p\in A_p)=E\left\{\prod_{i=1}^pP(X_i\in A_i\mid Z)\right\}=\int_\Theta\prod_{i=1}^pP_i(A_i\mid\theta)\,\pi(d\theta)$$
for all $A_1,\ldots,A_p\in\mathcal{B}(\mathbb{R})$. Define a probability measure $\lambda$ on $\mathbb{R}^{2p}$ as
$$\lambda(A_1\times\ldots\times A_p\times B_1\times\ldots\times B_p)=\int_\Theta\prod_{i=1}^pP_i(A_i\mid\theta)\,P_i(B_i\mid\theta)\,\pi(d\theta)$$
where $A_i,\,B_i\in\mathcal{B}(\mathbb{R})$ for each $i$. Then, as shown in \cite{JSPI}, any $p$-variate random vector $\widetilde{X}$ such that $\mathcal{L}(X,\widetilde{X})=\lambda$ is a knockoff copy of $X$. Such an $\widetilde{X}$ is called a {\em conditional independence knockoff} (CIK); see \cite{EJS}.

\medskip

\noindent To illustrate, suppose that $P_i(\cdot\mid\theta)$ admits a density $f_i(\cdot\mid\theta)$ with respect to some dominating measure $\gamma_i$. For instance, $\gamma_i$ could be Lebesgue measure for some $i$ and counting measure for some other $i$. Then, $X$ has a density with respect to the product measure
$$\gamma=\gamma_1\times\ldots\times\gamma_p.$$
Such a density can be written as
$$h(x)=h(x_1,\ldots,x_p)=\int_\Theta\prod_{i=1}^pf_i(x_i\mid\theta)\,\pi(d\theta)\quad\quad\text{for }x\in\mathbb{R}^p.$$
Similarly, $(X,\widetilde{X})$ has a density with respect to $\gamma\times\gamma$, that is,
\begin{gather*}
f(x,\widetilde{x})=f(x_1,\ldots,x_p,\widetilde{x}_1,\ldots,\widetilde{x}_p)=\int_\Theta\prod_{i=1}^pf_i(x_i\mid\theta)\,f_i(\widetilde{x}_i\mid\theta)\,\pi(d\theta)
\end{gather*}
where $x,\,\widetilde{x}\in\mathbb{R}^p$. Hence, the conditional density of $\widetilde{X}$ given $X=x$ is
$$\frac{f(x,\widetilde{x})}{h(x)}=\frac{\int_\Theta\prod_{i=1}^pf_i(x_i\mid\theta)\,f_i(\widetilde{x}_i\mid\theta)\,\pi(d\theta)}{\int_\Theta\prod_{i=1}^pf_i(x_i\mid\theta)\,\pi(d\theta)}\quad\text{whenever }h(x)>0.$$

\bigskip

\noindent The price to be paid, for a CIK $\widetilde{X}$ to be available, is condition \eqref{m9i}. In exchange of this price, there are at least three advantages. Firstly, to build $\widetilde{X}$ is straightforward. Secondly, $\mathcal{L}(X,\widetilde{X})$ can be written in closed form. Thirdly, $\widetilde{X}$ is sometimes optimal under some popular criterions. For instance, if $Z$ is such that $E(X_i\mid Z)=0$, one obtains
\begin{gather*}
\text{Cov}(X_i,\widetilde{X}_i)=E(X_i\widetilde{X}_i)-E(X_i)^2
\\=E\bigl\{E(X_i\mid Z)E(\widetilde{X}_i\mid Z)\bigr\}-E\bigl\{E(X_i\mid Z)\bigr\}^2=E\bigl\{E(X_i\mid Z)^2\bigr\}=0.
\end{gather*}
Hence, $\widetilde{X}$ is best possible if the optimality criterion is minimizing the mean absolute correlation. Similarly, as shown in \cite[Example 1]{EJS}, $\widetilde{X}$  is optimal if the criterion is minimizing the reconstructability index of \cite{SPJA}.

\medskip

\section{Exchangeable binary covariates}\label{b76t}

\noindent As discussed in Section \ref{intro}, under condition \eqref{i9p} (i.e., when the covariates are categorical), exchangeability of $X$ can be supported by various arguments. In this section, we actually assume condition \eqref{i9p} with $F=\{0,1\}$ and we take $X$ (infinitely) exchangeable. Equivalently, we assume
\begin{gather}\label{ni9x2}
P(X=x)=\int_{[0,1]}u^{n_1}(1-u)^{n_0}\,\pi(du),\quad\quad x\in\{0,1\}^p,
\end{gather}
for some (unique) probability measure $\pi$ on $[0,1]$. Such a $\pi$ is called ``the prior" in the sequel. Recall also that
$$n_j=n_j(x)=\sum_{i=1}^p\textbf{1}(x_i=j)\quad\quad\text{for all }j\in\{0,1\}\text{ and }x\in\{0,1\}^p.$$

\medskip

\noindent Under condition \eqref{ni9x2}, building a CIK is straightforward.

\begin{theorem}\label{gy78n3} Assume condition \eqref{ni9x2} and denote by $\widetilde{X}$ any $p$-variate random vector. Then, $\widetilde{X}$ is a knockoff provided
\begin{gather*}
P\bigl(X=x,\,\widetilde{X}=\widetilde{x}\bigr)=\int_{[0,1]}u^{n_1+\widetilde{n}_1}(1-u)^{n_0+\widetilde{n}_0}\,\pi(du)
\end{gather*}
for all $x,\,\widetilde{x}\in\{0,1\}^p$, where $\widetilde{n}_j=n_j(\widetilde{x})=\sum_{i=1}^p\textbf{1}(\widetilde{x}_i=j)$.
\end{theorem}

\begin{proof}
Let $Z,D_1,\ldots,D_p$ be random variables such that $Z\sim\pi$, $Z$ is independent of $(D_1,\ldots,D_p)$, and $D_1,\ldots,D_p$ are i.i.d. with uniform distribution on $[0,1]$. Define
$$X_i^*=\textbf{1}(D_i\le Z)\quad\quad\text{for }i=1,\ldots,p.$$
By condition \eqref{ni9x2}, $(X_1^*,\ldots,X_p^*)$ has the same distribution as $X$. Hence, in order to prove the theorem, it can be assumed $X_i=X_i^*$ for each $i$. In this case, $X_1,\ldots,X_p$ are conditionally i.i.d. given $Z$ with
$$P(X_1=1\mid Z)=P(D_1\le Z\mid Z)=Z.$$
Therefore, in the notation of Section \ref{vh87}, one can take $\Theta=[0,1]$ and $P_i(\{1\}\mid\theta)=\theta$ for all $i$.
\end{proof}

\medskip

\noindent Note that Theorem \ref{gy78n3} provides an explicit formula for the conditional distribution of $\widetilde{X}$ given $X=x$, that is
\begin{gather}\label{bg78n3w}
P(\widetilde{X}=\widetilde{x}\mid X=x)=\frac{P(X=x,\,\widetilde{X}=\widetilde{x})}{P(X=x)}=\frac{\int_{[0,1]}u^{n_1+\widetilde{n}_1}(1-u)^{n_0+\widetilde{n}_0}\,\pi(du)}{\int_{[0,1]}u^{n_1}(1-u)^{n_0}\,\pi(du)}
\end{gather}
whenever $P(X=x)>0$.

\medskip

\noindent If one agrees to assume $X$ exchangeable and to use the knockoff $\widetilde{X}$ of Theorem \ref{gy78n3}, the only remaining task is to select the prior $\pi$. We now consider a few specific choices of $\pi$. We stress that our goal is not to introduce new sophisticated priors (there is a huge literature on that) but to indicate some priors which make KP effective in real problems. We also note that the numerical performances of such priors are investigated in Section \ref{b5f9k}. Here, we just evaluate $\mathcal{L}(X)$ and $\mathcal{L}(X,\widetilde{X})$.

\medskip

\noindent We begin with discrete priors. Let
$$I=\left\{0,\frac{1}{p},\frac{2}{p},\ldots,1\right\}\quad\text{and}\quad\pi\{u\}=P\left(\sum_{i=1}^pX_i=pu\right)\quad\text{for each }u\in I.$$
This choice of $\pi$ can be supported in at least two ways. First, suppose the statistician selects an arbitrary prior $\pi^*$. Then, under mild conditions on $\pi^*$, there are two constants $a>0$ and $b>0$ such that
$$\frac{a}{p}\le d(\pi,\pi^*)\le\frac{b}{p}\,,$$
where $d$ is the bounded Lipschitz metric; see \cite{BPR17, MPS2016}. Thus, for large $p$, the probabilities $\pi$ and $\pi^*$ are close to each other (according to the distance $d$) and replacing $\pi^*$ with $\pi$ yields a negligible error. Second, and more important, $\pi$ admits a simple interpretation which provides a criterion to select it in real problems. In fact, $\pi\{u\}=P(A_u)$ where
$$A_u=\bigl\{\text{exactly }pu\text{ covariates are equal to }1\bigr\}.$$
Hence, to select $\pi$, the statistician is only required to assign $P(A_u)$ for each $u\in I$. In genetic applications, for instance, he/she has to evaluate the probability of observing exactly $pu$ genes among the $p$ possible ones. In addition to be reasonable, this request is strictly connected to the specific features of the problem at hand.

\medskip

\noindent In case of complete ignorance, $\pi$ can be taken to be uniform on $I$ or on a subset $J\subset I$. As an example, for large $p$, the two extreme cases $u=0$ and $u=1$ are often negligible. In this case, one can let $J=\bigl\{1/p,\ldots,(p-1)/p\bigr\}$ and $\pi\{u\}=1/\abs{J}$ for each $u\in J$. Another obvious choice of $\pi$ is
\begin{gather}\label{bin}
\pi\{u\}=\left(
             \begin{array}{c}
               p \\
               pu \\
             \end{array}
           \right)\,\alpha^{pu}\,(1-\alpha)^{p(1-u)}\quad\quad\text{for each }u\in I,
\end{gather}
where $\alpha\in (0,1)$ is a fixed number. Roughly speaking, $\alpha$ can be regarded as the probability that a generic covariate turns out to be 1 in a single trial. Once again, in case of ignorance, it is tempting to let $\alpha=1/2$. In any case, for each prior $\pi$ such that $\pi(I)=1$, one obtains
\begin{gather*}
P(X=x)=\sum_{u\in I}\pi\{u\}\,u^{n_1}\,(1-u)^{n_0}\quad\quad\text{and}
\\P\bigl(X=x,\,\widetilde{X}=\widetilde{x}\bigr)=\sum_{u\in I}\pi\{u\}\,u^{n_1+\widetilde{n}_1}\,(1-u)^{n_0+\widetilde{n}_0}.
\end{gather*}

\medskip

\noindent Discrete priors have a nice interpretation and are useful in real problems. However, when dealing with exchangeable sequences of indicators, it is quite usual to take $\pi$ absolutely continuous, i.e. $\pi(du)=f(u)\,du$ for some density $f$ on $[0,1]$. As noted above, for our purposes, there is no reason for an involved choice of $f$. Hence, we take $f$ to be a beta density with parameters $a>0$ and $b>0$. With this choice of $f$, one immediately obtains
\begin{gather*}
P(X=x)=\int_0^1u^{n_1}(1-u)^{n_0}\,f(u)\,du
=\frac{\Gamma(a+b)}{\Gamma(a)\,\Gamma(b)}\,\int_0^1 u^{a+n_1-1}(1-u)^{b+n_0-1}\,du
\\=\frac{\Gamma(a+b)\,\Gamma(n_0+b)\,\Gamma(n_1+a)}{\Gamma(a)\,\Gamma(b)\,\Gamma(p+a+b)}
\end{gather*}
and
\begin{gather*}
P\bigl(X=x,\,\widetilde{X}=\widetilde{x}\bigr)=\frac{\Gamma(a+b)\,\Gamma(n_0+\widetilde{n}_0+b)\,\Gamma(n_1+\widetilde{n}_1+a)}{\Gamma(a)\,\Gamma(b)\,\Gamma(2p+a+b)}.
\end{gather*}

\medskip

\noindent A last remark is in order.

\medskip

\begin{remark}\label{vg778nj} Recall that, for KP to be effective in terms of power, $X$ and $\widetilde{X}$ should be ``as independent as possible". For fixed $i$, since $X_i$ and $\widetilde{X}_i$ are indicators, $\text{Cov}(X_i,\widetilde{X}_i)$ provides a possible evaluation of how close $\mathcal{L}(X_i,\widetilde{X}_i)$ is to the independence law. Moreover, if $\widetilde{X}$ is the knockoff of Theorem \ref{gy78n3},
\begin{gather*}
\text{Cov}(X_i,\widetilde{X}_i)=E(X_i\widetilde{X}_i)-E(X_i)^2=\int_{[0,1]}u^2\,\pi(du)-\left(\int_{[0,1]}u\,\pi(du)\right)^2.
\end{gather*}
This suggests to select a prior $\pi$ with small variance. This is just a rough indication, however, for KP does not work in the extreme case where the variance of $\pi$ is null.
\end{remark}

\medskip

\section{Partially exchangeable binary covariates}\label{bv5r4}
\noindent In this section, we still assume $F=\{0,1\}$ but we take $X$ partially exchangeable. Precisely, there are an integer $k$ and a random vector $(U,V)$ such that

\medskip

\begin{itemize}

\item[(j)] $1\le k<p$, $0\le U\le 1$, $0\le V\le 1$;

\medskip

\item[(jj)] $X_1,\ldots,X_p$ are conditionally independent given $(U,V)$;

\medskip

\item[(jjj)] $P(X_i=1\mid U,V)=U$ if $i\le k$ and $P(X_i=1\mid U,V)=V$ if $i>k$.

\end{itemize}

\medskip

\noindent As noted in Section \ref{intro}, this modeling of $X$ may be reasonable when the available information leads to distinguish the subvectors $X'=(X_1,\ldots,X_k)$ and $X''=(X_{k+1},\ldots,X_p)$ and to regard $X'$ and $X''$ as marginally exchangeable. In this sense, our modeling may be related to \cite{KATSAB}, a main difference being that we require $P(X\in\{0,1\}^p)=1$.

\medskip

\noindent Let $\pi$ be the probability distribution of $(U,V)$. Define also
$$s_j=s_j(x)=\sum_{i=1}^k\textbf{1}(x_i=j)\quad\text{and}\quad t_j=t_j(x)=\sum_{i=k+1}^p\textbf{1}(x_i=j)$$
for all $x\in\{0,1\}^p$. In this notation, conditions (j)-(jj)-(jjj) imply
\begin{gather*}
P(X=x)=E\left\{\prod_{i=1}^pP(X_i=x_i\mid U,V)\right\}=E\bigl\{U^{s_1}(1-U)^{s_0}V^{t_1}(1-V)^{t_0}\bigr\}
\\=\int_{[0,1]^2}u^{s_1}(1-u)^{s_0}v^{t_1}(1-v)^{t_0}\,\pi(du,dv).
\end{gather*}
In addition, to get a CIK, it suffices to apply the next result.

\begin{theorem} Assume conditions (j)-(jj)-(jjj) and denote by $\pi$ the probability distribution of $(U,V)$ and by $\widetilde{X}$ any $p$-variate random vector. Then, $\widetilde{X}$ is a knockoff provided
\begin{gather*}
P\bigl(X=x,\,\widetilde{X}=\widetilde{x}\bigr)=\int_{[0,1]^2}u^{s_1+\widetilde{s}_1}(1-u)^{s_0+\widetilde{s}_0}\,v^{t_1+\widetilde{t}_1}(1-v)^{t_0+\widetilde{t}_0}\,\pi(du,dv)
\end{gather*}
for all $x,\,\widetilde{x}\in\{0,1\}^p$, where
$$\widetilde{s}_j=s_j(\widetilde{x})=\sum_{i=1}^k\textbf{1}(\widetilde{x}_i=j)\quad\text{and}\quad\widetilde{t}_j=t_j(\widetilde{x})=\sum_{i=k+1}^p\textbf{1}(\widetilde{x}_i=j).$$
\end{theorem}

\begin{proof}
Apply the theory of Section \ref{vh87} with $Z=(U,V)$ and $\Theta=[0,1]^2$. Moreover, for all $\theta=(u,v)\in [0,1]^2$, define
$$P_i(\{1\}\mid\theta)=u\,\,\text{ if }i\le k\quad\text{and}\quad P_i(\{1\}\mid\theta)=v\,\,\text{ if }i>k.$$
\end{proof}

\medskip

\noindent It remains to select $\pi$. As regards discrete priors, one option is choosing $\pi$ such that $\pi(I)=1$ where now $I$ is given by
$$I=\left\{\left(\frac{r}{k},\,\frac{s}{p-k}\right):r=0,1,\ldots,k;\,s=0,1,\ldots,p-k\right\}.$$
This type of priors are usually simple and admit a straightforward interpretation. If $(u,v)\in I$, in fact, $\pi\{(u,v)\}$ can be regarded as the probability of observing exactly $ku$ covariates equal to 1 within $X'$ and exactly $(p-k)v$ covariates equal to 1 within $X''$.

\medskip

\noindent Once again, in case of complete ignorance, $\pi$ can be taken to be uniform on $I$ or on a subset $J\subset I$. In analogy with \eqref{bin}, another option is to fix a number $\alpha\in (0,1)$, a Borel function $f:[0,1]\rightarrow [0,1]$, and to take $(U,V)$ such that
$$kU\sim\,\text{Bin}(k,\alpha)\quad\text{and}\quad (p-k)V\mid U\sim\,\text{Bin}(p-k,\,f(U)).$$
In this case, for each $(u,v)\in I$, the prior $\pi$ can be written as
\begin{gather*}
\pi\{(u,v)\}=P(U=u,\,V=v)=P(V=v\mid U=u)\,P(U=u)
\\=\left(
             \begin{array}{c}
               p-k \\
               (p-k)v \\
             \end{array}
           \right)\,f(u)^{(p-k)v}(1-f(u))^{(p-k)(1-v)}\left(
             \begin{array}{c}
               k \\
               ku \\
             \end{array}
           \right)\,\alpha^{ku}(1-\alpha)^{k(1-u)}.
\end{gather*}
To use this prior, one has to select $\alpha$ and $f$. Here, $\alpha$ should be regarded as the probability that a generic covariate is 1 in a single trial from $X'$. In turn, assuming that exactly $ku$ covariates of $X'$ are equal to 1, $f(u)$ should be viewed as the conditional probability that a generic covariate takes value 1 in a single trial from $X''$. In real problems, it is probably convenient to make simple choices of $f$ such as $f(u)=1-u$ or $f(u)=u^b$ for some $b>0$.

\medskip

\noindent In any case, if $\pi$ is such that $\pi(I)=1$, then
\begin{gather*}
P(X=x)=\sum_{(u,v)\in I}\pi\{(u,v)\}\,u^{s_1}\,(1-u)^{s_0}\,v^{t_1}\,(1-v)^{t_0}\quad\quad\text{and}
\\P\bigl(X=x,\,\widetilde{X}=\widetilde{x}\bigr)=\sum_{(u,v)\in I}\pi\{(u,v)\}\,u^{s_1+\widetilde{s}_1}\,(1-u)^{s_0+\widetilde{s}_0}\,v^{t_1+\widetilde{t}_1}\,(1-v)^{t_0+\widetilde{t}_0}.
\end{gather*}

\medskip

\noindent A probability measure which vanishes on singletons is said to be {\em diffuse}. Regarding diffuse priors, we just mention a straightforward choice of $\pi$ which allows simple calculations and may work for KP. Fix again a Borel function $f:[0,1]\rightarrow [0,1]$ and a diffuse probability measure $\mu$ on $[0,1]$. After selecting $f$ and $\mu$, take $\pi$ supported by the graph of $f$ and having marginal $\mu$ in the first coordinate. Equivalently, take $U\sim\mu$ and $V=f(U)$. In this case,
$$P(X=x)=\int_{[0,1]}u^{s_1}(1-u)^{s_0}f(u)^{t_1}(1-f(u))^{t_0}\,\mu(du).$$
As to $f$, one can choose a second probability $\mu^*$ on $[0,1]$ and take $f$ such that $f(U)\sim\mu^*$. Or else, as above, one can make simple choices such as $f(u)=1-u$ or $f(u)=u^b$ for some $b>0$. Note that, in the latter case, $f(u)\le u$ or $f(u)\ge u$ according to whether $b\ge 1$ or $b\le 1$. Hence, through $b$, the statistician can tune the probability of observing a covariate equal to 1 within $X'$ and within $X''$. Moreover, the evaluation of $\mathcal{L}(X,\,\widetilde{X})$ is straghtforward. For instance, if $f(u)=1-u$ and $\mu$ is the uniform distribution on $[0,1]$, then
\begin{gather*}
P(X=x)=\frac{\Gamma(s_1+t_0+1)\,\Gamma(s_0+t_1+1)}{\Gamma(p+2)}\quad\quad\text{and}
\\P\bigl(X=x,\,\widetilde{X}=\widetilde{x}\bigr)=\frac{\Gamma(s_1+\widetilde{s}_1+t_0+\widetilde{t}_0+1)\,\Gamma(s_0+\widetilde{s}_0+t_1+\widetilde{t}_1+1)}{\Gamma(2p+2)}.
\end{gather*}

\medskip

\noindent Partial exchangeability makes sense in a number of situations. We close this section with an example where condition \eqref{i9p} is satisfied by some, but not all, the covariates. This situation is actually quite common in applications; see e.g. \cite{STATMED}.

\medskip

\begin{example}\label{g6yh8n}  Given an integer $1\le k<p$, suppose
$$X_i\in\{0,1\}\quad\text{for }i\le k\quad\text{while}\quad X_i\text{ has a diffuse distribution for }i>k.$$
To model this situation, we first choose a collection $\{Q_v:v\in\mathcal{V}\}$ where $\mathcal{V}$ in any index set and $Q_v$ a diffuse probability measure on $\mathbb{R}$. Then, we take $X_1,\ldots,X_p$ conditionally independent, given a random vector $(U,V)$, with
$$P(X_i=1\mid U,V)=U\text{ if }i\le k\quad\text{and}\quad P(X_i\in\cdot\mid U,V)=Q_V(\cdot)\text{ if }i>k.$$
Then, for all $x_1,\ldots,x_k\in\{0,1\}$ and $A_{k+1},\ldots,A_p\in\mathcal{B}(\mathbb{R})$, one obtains
\begin{gather*}
P\bigl(X_i=x_i\text{ and }X_j\in A_j\text{ for }i=1,\ldots,k\text{ and }j=k+1,\ldots,p\bigr)
\\=\int_{[0,1]\times\mathcal{V}}u^{s_1}(1-u)^{s_0}\,\prod_{j=k+1}^pQ_v(A_j)\,\pi(du,dv)
\end{gather*}
where $\pi$ is the probability distribution of $(U,V)$ and $s_j=\sum_{i=1}^k\textbf{1}(x_i=j)$. Moreover, a CIK is any $p$-variate random vector $\widetilde{X}$ satisfying
\begin{gather*}
P\bigl(X_i=x_i,\,\widetilde{X}_i=\widetilde{x}_i,\,X_j\in A_j,\,\widetilde{X}_j\in B_j\text{ for }i=1,\ldots,k\text{ and }j=k+1,\ldots,p\bigr)
\\=\int_{[0,1]\times\mathcal{V}}u^{s_1+\widetilde{s}_1}(1-u)^{s_0+\widetilde{s}_0}\,\prod_{j=k+1}^pQ_v(A_j)\,Q_v(B_j)\,\pi(du,dv)
\end{gather*}
where $\widetilde{s}_j=\sum_{i=1}^k\textbf{1}(\widetilde{x}_i=j)$.

\medskip

\noindent To be more concrete, suppose
$$\mathcal{V}=\mathbb{R}\times (0,\infty),\quad v=(v_1,v_2)\quad\text{and}\quad Q_v=N(v_1,v_2).$$
Then, $X$ has a density $h$ with respect to $\gamma$, where
$$\gamma=\Bigl(\text{counting measure on }\{0,1\}^k\Bigr)\,\times\,\Bigl(\text{Lebesgue measure on }\mathbb{R}^{p-k}\Bigr).$$
Precisely, if $\phi(\cdot\mid v)$ denotes the density of $N(v_1,v_2)$, one obtains
\begin{gather}\label{9n6yh8k}
h(x)=\int_{[0,1]\times\mathcal{V}}u^{s_1}(1-u)^{s_0}\,\prod_{j=k+1}^p\phi(x_j\mid v)\,\pi(du,dv)
\end{gather}
where $x\in\{0,1\}^k\times\mathbb{R}^{p-k}$; see Section \ref{vh87}. For instance, if $U$ is uniform on $[0,1]$ and $V=(f_1(U),f_2(U))$, for some functions $f_1$ and $f_2$, equation \eqref{9n6yh8k} reduces to
$$h(x)=\int_0^1 u^{s_1}(1-u)^{s_0}\,\prod_{j=k+1}^p\phi(x_j\mid f_1(u),f_2(u))\,du.$$
\end{example}

\medskip

\section{Exchangeable categorical covariates}\label{b5rtg8n9}
\noindent We now turn to the general case $F=\bigl\{0,1,\ldots,m\bigr\}$. As an obvious example, in several genetic applications, it is reasonable to let $m=2$ and $X_i=0$, $X_i=1$ or $X_i=2$ according to whether the $i$-th gene is absent, present or unresponsive. We take $X$ infinitely exchangeable, but, at the price of minor complications, exchangeability could be weakened into partial exchangeability.

\medskip

\noindent This time, we assume that there is a random vector
$$U=(U_1,\ldots,U_m)$$
such that

\medskip

\begin{itemize}

\item[(a)] $U\in S$ where $S=\bigl\{u\in\mathbb{R}^m:u_j\ge 0\text{ for all }j\text{ and }\sum_{j=1}^mu_j\le 1\bigr\}$;

\medskip

\item[(b)] $X_1,\ldots,X_p$ are conditionally i.i.d. given $U$;

\medskip

\item[(c)] $P(X_1=j\mid U)=U_j\,$ for all $j\in F\setminus\{0\}$.
\end{itemize}

\medskip

\noindent It follows that, for each $x\in F^p$,
\begin{gather*}
P(X=x)=E\left\{\prod_{i=1}^pP(X_i=x_i\mid U)\right\}=E\left\{\Bigl(1-\sum_{j=1}^mU_j\Bigr)^{n_0}\,\prod_{j=1}^mU_j^{n_j}\right\}
\\=\int_S\,\,\Bigl(1-\sum_{j=1}^mu_j\Bigr)^{n_0}\,\prod_{j=1}^mu_j^{n_j}\,\pi(du)
\end{gather*}
where $\pi$ is the probability distribution of $U$. Moreover, a CIK is given by the following result.

\begin{theorem}\label{n89m3ew} Assume conditions (a)-(b)-(c) and denote by $\pi$ the probability distribution of $U$ and by $\widetilde{X}$ any $p$-variate random vector. Then, $\widetilde{X}$ is a knockoff provided
\begin{gather*}
P(X=x,\,\widetilde{X}=\widetilde{x}\bigr)=\int_S\,\,\Bigl(1-\sum_{j=1}^mu_j\Bigr)^{n_0+\widetilde{n}_0}\,\prod_{j=1}^mu_j^{n_j+\widetilde{n}_j}\,\pi(du)
\end{gather*}
for all $x,\,\widetilde{x}\in F^p$, where $\widetilde{n}_j=n_j(\widetilde{x})=\sum_{i=1}^p\textbf{1}(\widetilde{x}_i=j)$.
\end{theorem}

\begin{proof}
Obviously, $P(X_1=0\mid U)=1-\sum_{j=1}^mU_j$. Hence, once again, we can apply the theory of Section \ref{vh87} with $Z=U$ and $\Theta=S$. Precisely, for all $i=1,\ldots,p$ and all $\theta=(u_1,\ldots,u_m)\in S$, we let
$$P_i(\{0\}\mid\theta)=1-\sum_{j=1}^mu_j\quad\text{and}\quad P_i(\{j\}\mid\theta)=u_j\quad\text{if }j\in F\setminus\{0\}.$$
\end{proof}

\medskip

\noindent As to the choice of $\pi$, it is still convenient to distinguish discrete and diffuse priors. The former are meant as those priors $\pi$ such that $\pi(I)=1$ where
$$I=\left\{\left(\frac{r_1}{p},\ldots,\frac{r_m}{p}\right):r_1,\ldots,r_m\text{ non-negative integers and }\sum_{j=1}^mr_j\le p\right\}.$$
The interpretation of discrete priors is still the same, namely, for $u\in I$, $\pi\{u\}$ should be regarded as the probability of observing exactly $pu_j$ covariates equal to $j$ for each $j\in F\setminus\{0\}$. This interpretation provides an useful criterion for selecting $\pi$ in real problems. We do not write explicit formulae for not to be boring, but they should be clear at this point.

\medskip

\noindent Finally, regarding diffuse priors, we still make a standard choice and we take $\pi$ to be the Dirichlet distribution on $S$ with parameters $\alpha_0,\alpha_1,\ldots,\alpha_m>0$. This means that $\pi$ is the probability distribution of the vector
$$\left(\frac{G_1}{\sum_{j=0}^mG_j}\,,\,\ldots\,,\,\frac{G_m}{\sum_{j=0}^mG_j}\right)$$
where $G_0,G_1,\ldots,G_m$ are independent Gamma random variables and each $G_j$ has scale parameter 1 and shape parameter $\alpha_j$. With this choice of $\pi$, letting $\alpha=\sum_{j=0}^m\alpha_j$, one obtains
\begin{gather*}
P(X=x)=\frac{\Gamma(\alpha)}{\Gamma(\alpha_0)\ldots\Gamma(\alpha_m)}\,\,\int_S \Bigl(1-\sum_{j=1}^mu_j\Bigr)^{n_0+\alpha_0-1}\prod_{j=1}^mu_j^{n_j+\alpha_j-1}\,du
\\=\frac{\Gamma(\alpha)}{\Gamma(\alpha_0)\ldots\Gamma(\alpha_m)\,\Gamma(p+\alpha)}\,\,\prod_{j=0}^m\Gamma(n_j+\alpha_j)
\end{gather*}
and
$$P(X=x,\,\widetilde{X}=\widetilde{x}\bigr)=\frac{\Gamma(\alpha)}{\Gamma(\alpha_0)\ldots\Gamma(\alpha_m)\,\Gamma(2p+\alpha)}\,\,\prod_{j=0}^m\Gamma(n_j+\widetilde{n}_j+\alpha_j).$$

\medskip

\section{Robustness of the CIKs with respect to the prior}\label{rob1}
\noindent Suppose one agrees to assume $X$ infinitely exchangeable and aims to calculate a CIK $\widetilde{X}$ (that is, to evaluate a knockoff $\widetilde{X}$ via Theorem \ref{n89m3ew}). To do this, a prior $\pi$ is to be selected and the subsequent statistical analysis depends on $\pi$. A quite natural question is: ``how much this analysis depends on the choice of $\pi$ ?" Giving a {\em quantitative} answer is not easy. Hence, we try to answer a related question, namely: ``how much $\mathcal{L}(X,\widetilde{X})$ and $\mathcal{L}(\widetilde{X}\mid X=x)$ depend on $\pi$ ?"

\medskip

\noindent In this section, we let $F=\{0,1\}$ and we assume $X$ infinitely exchangeable. Accordingly, a prior is a probability measure on $[0,1]$. We also recall that, if $P_1$ and $P_2$ are probability measures on any measurable space $(\mathcal{X},\mathcal{B})$, their {\em total variation distance} is
$$\norm{P_1-P_2}=\sup_{B\in\mathcal{B}}\,\abs{P_1(B)-P_2(B)}.$$

\medskip

\noindent Our starting point is the following result.

\begin{lemma}\label{i9m0}
Let $T=(T_1,\ldots,T_n)$ be a $n$-variate random vector of indicators (i.e., $T_i\in\{0,1\}$ for all $i=1,\ldots,n$). If $T$ is infinitely exchangeable, then
$$\norm{\mathcal{L}_{\pi_1}(T)-\mathcal{L}_{\pi_2}(T)}\le\norm{\pi_1-\pi_2}$$
for any priors $\pi_1$ and $\pi_2$, where $\mathcal{L}_\pi(T)$ denotes the probability distribution of $T$ when the prior is $\pi$.
\end{lemma}

\begin{proof}
We first recall that, if $P_1$ and $P_2$ are probability laws on a measurable space $(\mathcal{X},\mathcal{B})$, then $\norm{P_1-P_2}$ can be written as
$$\norm{P_1-P_2}=\sup_{0\le f\le 1,\,f\text{ measurable}}\,\Abs{\int_\mathcal{X} f\,dP_1-\int_\mathcal{X} f\,dP_2};$$
see e.g. \cite[Section 2]{PR24}. Having noted this fact, define $n_1(a)=\sum_{i=1}^na_i$ and $n_0(a)=n-n_1(a)$ for each $a\in F^n$. Define also
$$f_A(u)=\sum_{a\in A}u^{n_1(a)}(1-u)^{n_0(a)}\quad\text{for all }A\subset F^n\text{ and }u\in [0,1]$$
where the indeterminate form $0^0$ is meant as 1. Then, $f_A$  is a measurable function on $[0,1]$ for fixed $A$. In addition,
$$0\le f_A(u)\le\sum_{a\in F^n}u^{n_1(a)}(1-u)^{n_0(a)}=1.$$
It follows that
\begin{gather*}
\norm{\mathcal{L}_{\pi_1}(T)-\mathcal{L}_{\pi_2}(T)}=\sup_{A\subset F^n}\,\Abs{P_{\pi_1}\bigl(T\in A)-P_{\pi_2}(T\in A)}
\\=\sup_{A\subset F^n}\,\Abs{\,\sum_{a\in A}\,\int_{[0,1]} u^{n_1(a)}(1-u)^{n_0(a)}\,(\pi_1(du)-\pi_2(du))}
\\=\sup_{A\subset F^n}\,\Abs{\int_{[0,1]} f_A\,d\pi_1-\int_{[0,1]} f_A\,d\pi_2}
\\\le \sup_{0\le f\le 1, f\text{ measurable}}\,\Abs{\int_{[0,1]} f\,d\pi_1-\int_{[0,1]} f\,d\pi_2}=\norm{\pi_1-\pi_2}.
\end{gather*}
\end{proof}

\medskip

\noindent Letting $n=2p$ and $T=(X,\widetilde{X})$, Lemma \ref{i9m0} yields
\begin{gather*}
\norm{\mathcal{L}_{\pi_1}(X,\widetilde{X})-\mathcal{L}_{\pi_2}(X,\widetilde{X})}\le\norm{\pi_1-\pi_2}.
\end{gather*}

\medskip

\noindent Even if rough, the previous inequality provides a quantitative indication about the dependence of $\mathcal{L}_{\pi}(X,\widetilde{X})$ on the prior $\pi$. As regards KP, however, a more interesting object is $\mathcal{L}_{\pi}(\widetilde{X}\mid X=x)$, that is, the conditional distribution of $\widetilde{X}$ given $X=x$ when the prior is $\pi$. In fact, after observing $X=x$, the values $\widetilde{x}$ of $\widetilde{X}$ (required for applying KP) should be sampled from $\mathcal{L}_{\pi}(\widetilde{X}\mid X=x)$.

\medskip

\noindent To investigate $\mathcal{L}_{\pi}(\widetilde{X}\mid X=x)$, a first hint comes from Lemma \ref{i9m0} again. In fact, if $P_\pi(X=x)>0$, formula \eqref{bg78n3w} implies
\begin{gather*}
P_\pi(\widetilde{X}=\widetilde{x}\mid X=x)=\int_{[0,1]}u^{n_1(\widetilde{x})}(1-u)^{n_0(\widetilde{x})}\,\pi(du\mid x)
\\\text{where}\quad\pi(du\mid x)=\frac{u^{n_1(x)}(1-u)^{n_0(x)}}{P_\pi(X=x)}\,\pi(du).
\end{gather*}
Hence, $\mathcal{L}_{\pi}(\widetilde{X}\mid X=x)$ is an exchangeable law on $F^p$ with prior $\pi(\cdot\mid x)$. By Lemma \ref{i9m0}, applied with $n=p$ and $T=\widetilde{X}$, one obtains
$$\norm{\mathcal{L}_{\pi_1}(\widetilde{X}\mid X=x)-\mathcal{L}_{\pi_2}(\widetilde{X}\mid X=x)}\le\norm{\pi_1(\cdot\mid x)-\pi_2(\cdot\mid x)}.$$
Another upper bound, possibly more useful, is provided by the next result.

\begin{theorem}\label{v5c34g8}
If $X$ is infinitely exchangeable, then
\begin{gather}\label{x556b8uj9}
\norm{\mathcal{L}_{\pi_1}(\widetilde{X}\mid X=x)-\mathcal{L}_{\pi_2}(\widetilde{X}\mid X=x)}\le
\\\le\frac{\norm{\pi_1-\pi_2}+\Abs{P_{\pi_1}(X=x)-P_{\pi_2}(X=x)}}{\max\bigl\{P_{\pi_1}(X=x),\,P_{\pi_2}(X=x)\bigr\}}\notag
\end{gather}
for all $x\in F^p$ and all priors $\pi_i$ such that $P_{\pi_i}(X=x)>0$, $i=1,2$. Moreover,
$$\norm{\mathcal{L}_{\pi_1}(\widetilde{X}\mid X=x)-\mathcal{L}_{\pi_2}(\widetilde{X}\mid X=x)}\le\frac{2\,\norm{\pi_1-\pi_2}}{\max\bigl\{P_{\pi_1}(X=x),\,P_{\pi_2}(X=x)\bigr\}}.$$
\end{theorem}

\begin{proof}
First note that
$$\Abs{P_{\pi_1}(X=x)-P_{\pi_2}(X=x)}\le\norm{\mathcal{L}_{\pi_1}(X)-\mathcal{L}_{\pi_2}(X)}\le\norm{\pi_1-\pi_2}$$
where the second inequality depends on Lemma \ref{i9m0}. Hence, under \eqref{x556b8uj9}, one obtains
$$\norm{\mathcal{L}_{\pi_1}(\widetilde{X}\mid X=x)-\mathcal{L}_{\pi_2}(\widetilde{X}\mid X=x)}\le\frac{2\,\norm{\pi_1-\pi_2}}{\max\bigl\{P_{\pi_1}(X=x),\,P_{\pi_2}(X=x)\bigr\}}.$$
We next prove inequality \eqref{x556b8uj9}. For each $A\subset F^p$,
\begin{gather*}
\Abs{\,P_{\pi_1}(\widetilde{X}\in A\mid X=x)-P_{\pi_2}(\widetilde{X}\in A\mid X=x)}\le
\\\le\frac{\Abs{\,P_{\pi_1}(X=x,\,\widetilde{X}\in A)-P_{\pi_2}(X=x,\,\widetilde{X}\in A)}}{P_{\pi_1}(X=x)}\,+\,P_{\pi_2}(X=x,\,\widetilde{X}\in A)\,\frac{\Abs{P_{\pi_1}(X=x)-P_{\pi_2}(X=x)}}{P_{\pi_1}(X=x)\,P_{\pi_2}(X=x)}
\\\le\frac{\norm{\mathcal{L}_{\pi_1}(X,\widetilde{X})-\mathcal{L}_{\pi_2}(X,\widetilde{X})}}{P_{\pi_1}(X=x)}\,+\,\frac{\Abs{P_{\pi_1}(X=x)-P_{\pi_2}(X=x)}}{P_{\pi_1}(X=x)}
\\\le\frac{\norm{\pi_1-\pi_2}}{P_{\pi_1}(X=x)}\,+\,\frac{\Abs{P_{\pi_1}(X=x)-P_{\pi_2}(X=x)}}{P_{\pi_1}(X=x)}
\end{gather*}
where the last inequality still depends on Lemma \ref{i9m0}. By exactly the same argument,
\begin{gather*}
\Abs{\,P_{\pi_1}(\widetilde{X}\in A\mid X=x)-P_{\pi_2}(\widetilde{X}\in A\mid X=x)}\le\frac{\norm{\pi_1-\pi_2}}{P_{\pi_2}(X=x)}\,+\,\frac{\Abs{P_{\pi_1}(X=x)-P_{\pi_2}(X=x)}}{P_{\pi_2}(X=x)}.
\end{gather*}
Therefore,
\begin{gather*}
\norm{\mathcal{L}_{\pi_1}(\widetilde{X}\mid X=x)-\mathcal{L}_{\pi_2}(\widetilde{X}\mid X=x)}=\sup_{A\subset F^p}\,\Abs{\,P_{\pi_1}(\widetilde{X}\in A\mid X=x)-P_{\pi_2}(\widetilde{X}\in A\mid X=x)}
\\\le\frac{\norm{\pi_1-\pi_2}+\Abs{P_{\pi_1}(X=x)-P_{\pi_2}(X=x)}}{\max\bigl\{P_{\pi_1}(X=x),\,P_{\pi_2}(X=x)\bigr\}}.
\end{gather*}
This proves inequality \eqref{x556b8uj9} and concludes the proof of the theorem.
\end{proof}

\medskip

\noindent To exploit Theorem \ref{v5c34g8}, it is useful to have an estimate of $\norm{\pi_1-\pi_2}$. Such estimate is actually available in various meaningful cases. For instance, in view of Proposition 6 of \cite{KELB}, $\norm{\pi_1-\pi_2}$ can be upper bounded when $\pi_1$ and $\pi_2$ are both of the binomial type (i.e., they meet equation \eqref{bin}). Or else, when $\pi_1$ and $\pi_2$ are beta distributions. We close this section with an upper bound for $\norm{\pi_1-\pi_2}$ in the latter case. To this end, we denote by $\pi(a,b)$ the beta distribution with parameters $a$ and $b$.

\medskip

\begin{lemma}
Let $0<c<d$. If $a_j,\,b_j\in [c,d]$ for $j=1,2$, there is a constant $q$ depending only on $c$ and $d$ such that
$$\norm{\pi(a_1,b_1)-\pi(a_2,b_2)}\le q\,\bigl\{\abs{a_1-a_2}+\abs{b_1-b_2}\bigr\}.$$
\end{lemma}

\begin{proof}
It suffices to show that
\begin{gather}\label{ma6b9}
\norm{\pi(a_1,b)-\pi(a_2,b)}\le q\,\abs{a_1-a_2}\quad\quad\text{for all }a_1,\,a_2,\,b\in [c,d].
\end{gather}
In fact, under \eqref{ma6b9}, one obtains
\begin{gather*}
\norm{\pi(a_1,b_1)-\pi(a_2,b_2)}\le\norm{\pi(a_1,b_1)-\pi(a_2,b_1)}+\norm{\pi(a_2,b_1)-\pi(a_2,b_2)}
\\=\norm{\pi(a_1,b_1)-\pi(a_2,b_1)}+\norm{\pi(b_1,a_2)-\pi(b_2,a_2)}\le q\,\bigl\{\abs{a_1-a_2}+\abs{b_1-b_2}\bigr\}.
\end{gather*}
To prove \eqref{ma6b9}, we fix $a_1,\,a_2,\,b\in [c,d]$. Without loss of generality, we assume $a_1\le a_2$. Moreover, for all $u,\,v>0$, we let
$$\gamma(u,v)=\frac{\Gamma(u)\Gamma(v)}{\Gamma(u+v)}\quad\text{and}\quad\psi(u)=\frac{d}{du}\log \Gamma(u)=\frac{\Gamma'(u)}{\Gamma(u)}.$$
Then, for each $a>0$,
\begin{gather*}
\int_0^1 (1-u)^{b-1}\,u^{a-1}\,\log u\,du=\int_0^1 (1-u)^{b-1}\,\left(\frac{d}{da}u^{a-1}\right)\,du
\\=\frac{d}{da}\,\int_0^1 (1-u)^{b-1}\,u^{a-1}\,du=\frac{d}{da}\,\gamma(a,b)=-\gamma(a,b)\,\bigl\{\psi(a+b)-\psi(a)\bigr\}.
\end{gather*}
Since $1-u^{a_2-a_1}\le (a_1-a_2)\log u$ for all $u\in (0,1)$, it follows that
\begin{gather*}
\gamma(a_1,b)-\gamma(a_2,b)=\int_0^1 (1-u)^{b-1}\,u^{a_1-1}\,\bigl(1-u^{a_2-a_1}\bigr)\,du
\\\le (a_1-a_2)\,\int_0^1 (1-u)^{b-1}\,u^{a_1-1}\,\log u\,du
\\= (a_2-a_1)\,\gamma(a_1,b)\,\bigl\{\psi(a_1+b)-\psi(a_1)\bigr\}.
\end{gather*}
Finally, since $c\le a_1<a_1+b\le 2d$ and $\psi$ is strictly increasing, one obtains
\begin{gather*}
\norm{\pi(a_1,b)-\pi(a_2,b)}=\frac{1}{2}\,\int_0^1\Abs{\frac{u^{a_1-1}(1-u)^{b-1}}{\gamma(a_1,b)}-\frac{u^{a_2-1}(1-u)^{b-1}}{\gamma(a_2,b)}}\,du
\\\le\frac{1}{2\,\gamma(a_1,b)}\,\int_0^1 (1-u)^{b-1}\,u^{a_1-1}\,\bigl(1-u^{a_2-a_1}\bigr)\,du+\frac{\gamma(a_1,b)-\gamma(a_2,b)}{2\,\gamma(a_1,b)}
\\=\frac{\gamma(a_1,b)-\gamma(a_2,b)}{\gamma(a_1,b)}\le (a_2-a_1)\,\bigl\{\psi(2d)-\psi(c)\bigr\}.
\end{gather*}
Hence, \eqref{ma6b9} holds with $q=\psi(2d)-\psi(c)$.
\end{proof}

\medskip

\section{Numerical illustrations}\label{b5f9k}
\noindent This section reports a numerical experiment based on simulated and real data. The knockoffs obtained so far (henceforth, the CIKs) are compared with the {\em Model-$X$ Gaussian knockoff} of \cite{CFJL18} (henceforth, the Model-$X$). Indeed, the Model-$X$ is probably the gold standard among those knockoffs obtained by modeling $X$ as an absolutely continuous random vector, and this is why it has been taken as a reference. The comparison relies on the power and the false discovery rate. A further popular criterion would be the (observed) correlation between $X_i$ and $\widetilde{X}_i$. However, this criterion has been neglected. In fact, since the points of $F$ are merely labels, $\text{Cov}(X_i,\widetilde{X}_i)$ does not make sense in general. Possibly, $\text{Cov}(X_i,\widetilde{X}_i)$ could be informative when $F=\{0,1\}$, but in this case we already know that $\text{Cov}(X_i,\widetilde{X}_i)$ agrees with the variance of the prior $\pi$; see Remark \ref{vg778nj}.

\medskip

\subsection{Outline of the experiment}\label{pop98} In case of simulated data, we let $p=100$ and we focus on $F=\{0,1\}$ and $F=\{0,1,2\}$. The simulation experiment has been performed by adapting the standard format (see e.g. \cite{CFJL18, RSC}) to categorical covariates. Precisely:

\medskip

\begin{itemize}

\item A subset $T\subset\bigl\{1,\ldots,p\bigr\}$ such that $|T|=60$ has been randomly selected and the coefficients $\beta_1,\ldots,\beta_p$ have been defined as
$$\beta_i=0\,\text{ if }i\notin T\quad\text{and}\quad\beta_i=\frac{u}{\sqrt{n}}\,\text{ if }i\in T.$$
Here, $n$ is a positive integer and $u>0$ a parameter called {\em signal amplitude}.

\medskip

\item $n$ i.i.d. observations
$$X^{(j)}=\bigl(X_{1j},\ldots,X_{pj}\bigr),\quad\quad j=1,\ldots,n,$$
have been sampled from $\mathcal{L}(X)$.

\medskip

\item For each $j=1,\ldots,n$, the response variable $Y^{(j)}$ has been defined as follows. Let $e_1,\ldots,e_n$ be i.i.d. standard normal errors. If $F=\{0,1,2\}$, then
$$Y^{(j)}=\sum_{i=1}^p\Bigl\{2 \, \beta_i\,\textbf{1}(X_{ij}=0)+\beta_i \,\textbf{1}(X_{ij}=1) + \frac{\beta_i}{2} \,\textbf{1}(X_{ij}=2)\Bigr\}+e_j.$$
If $F=\{0,1\}$, the same model applies. However, $Y^{(j)}$ reduces to
$$Y^{(j)}=\sum_{i=1}^p\Bigl\{2 \, \beta_i\,\textbf{1}(X_{ij}=0)+\beta_i \,\textbf{1}(X_{ij}=1)\Bigr\}+e_j$$
since $\textbf{1}(X_{ij}=2)=0$ for all $i$ and $j$.

\medskip

\item For each $j=1,\ldots,n$, we sampled $m$ CIKs, say $\widetilde{X}^{(1,j)},\ldots,\widetilde{X}^{(m,j)}$, from the conditional distribution $\mathcal{L}(\widetilde{X}\mid X=x^{(j)})$ where $x^{(j)}$ is the observed value of $X^{(j)}$.

\item For each $k=1,\ldots,m$, KP has been applied to the data
$$\bigl\{Y^{(j)},\,X^{(j)},\,\widetilde{X}^{(k,j)}:j=1,\ldots,n\bigr\}$$
so as to calculate the power and the false discovery rate, say $POW(k)$ and $FDR(k)$. To do this, we exploited the R-cran package \verb"knockoff":

\medskip

\verb"https://cran.r-project.org/web/packages/knockoff/index.html".

\medskip

\noindent This package is based on the comparison between the lasso coefficient estimates of each covariate and its knockoff.

\medskip

\item The final outputs are the arithmetic means of the powers and the false discovery rates, i.e.,
$$POW=(1/m)\,\sum_{k=1}^m POW(k)\quad\text{and}\quad FDR=(1/m)\,\sum_{k=1}^m FDR(k).$$

\end{itemize}

\medskip

\noindent To run the simulation experiment, we took $m=n=1000$ and a theoretical value of the false discovery rate equal to $0.1$.

\medskip

\noindent In case of real data, the experiment is essentially the same. The main differences are that:

\medskip

\begin{itemize}

\item The dataset $\bigl\{Y^{(j)},\,X^{(j)}:j=1,\ldots,n\bigr\}$ is observed and not simulated;

\item $m=1$, i.e., exactly one knockoff $\widetilde{X}^{(j)}=\widetilde{X}^{(1,j)}$ is sampled for each $j$;

\item The theoretical value of the false discovery rate is $0.2$.

\end{itemize}

\medskip

\noindent The real application is based on the so called {\em human immunodeficiency virus type 1} (HIV-1) dataset, which has been already used in several previous studies on KP; see e.g. \cite{BC15, HIV, RSC}. The dimensions of HIV-1 are $n=846$ and $p=341$. The knockoff filter is applied to identify mutations associated with drug resistance. The HIV-1 includes drug resistance measurements and genotype information from HIV-1 samples, with separate datasets for resistance to protease inhibitors, nucleoside reverse transcriptase inhibitors, and non-nucleoside RT inhibitors. We focus on resistance to protease inhibitors and analyze the following drugs separately: amprenavir (APV), atazanavir (ATV), indinavir (IDV), lopinavir (LPV), nelfinavir (NFV), ritonavir (RTV), and saquinavir (SQV).

\medskip

\subsection{Exchangeable covariates}\label{x4r55b7} The performances of the CIKs and Model-$X$ in the numerical experiment are summarized by means of graphs. In each graph, the false discovery rate (FDR) is on the left and the power (POW) on the right. The CIKs correspond to a solid line while the Model-$X$ to a dashed line. The signal amplitude ranges in the set $\{1.5,2,3,4,5,7.5,10,15,20\}$.

\medskip

\begin{figure}[!htbp]
        \centering
        \hbox{
        \hspace{-2cm}\includegraphics[width=0.7\textwidth]{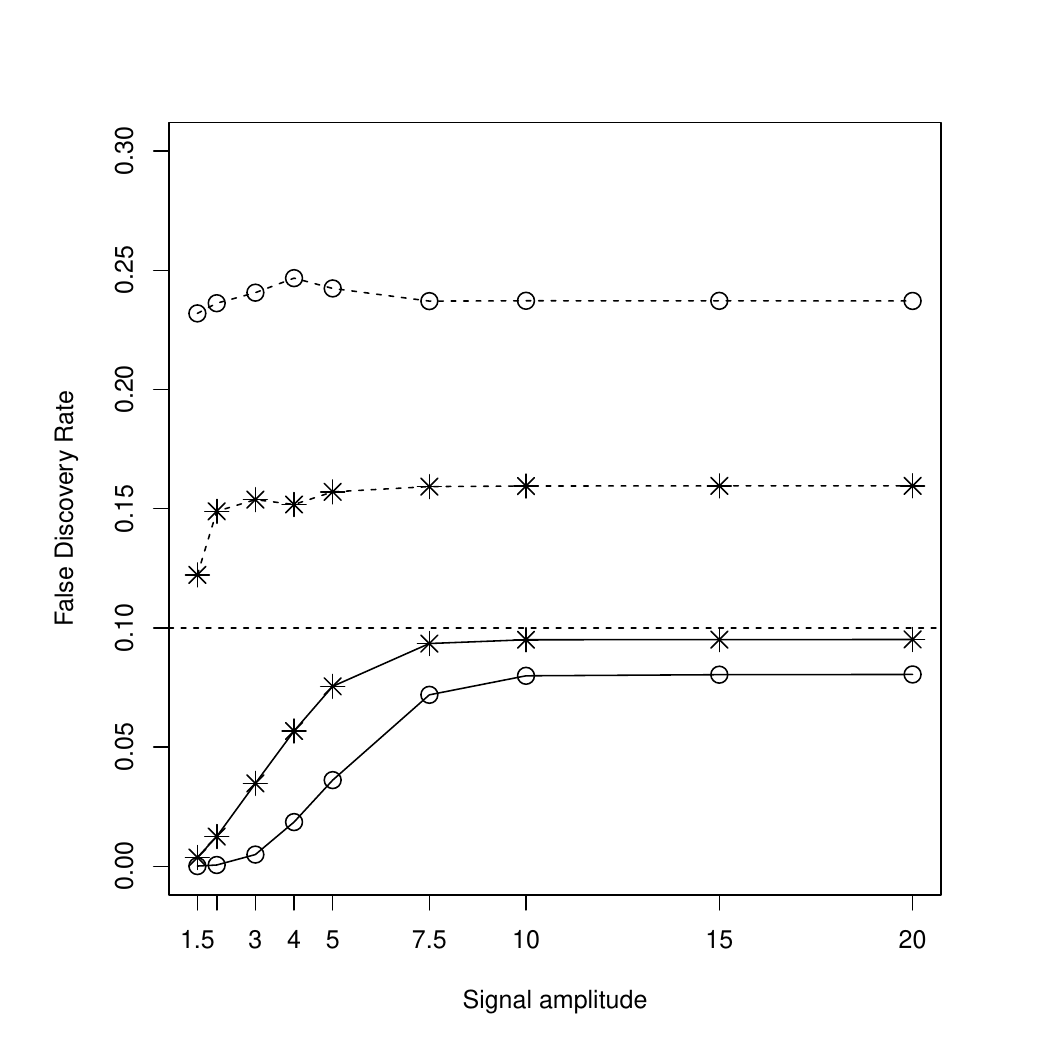}\hspace{-1cm}
         \includegraphics[width=0.7\textwidth]{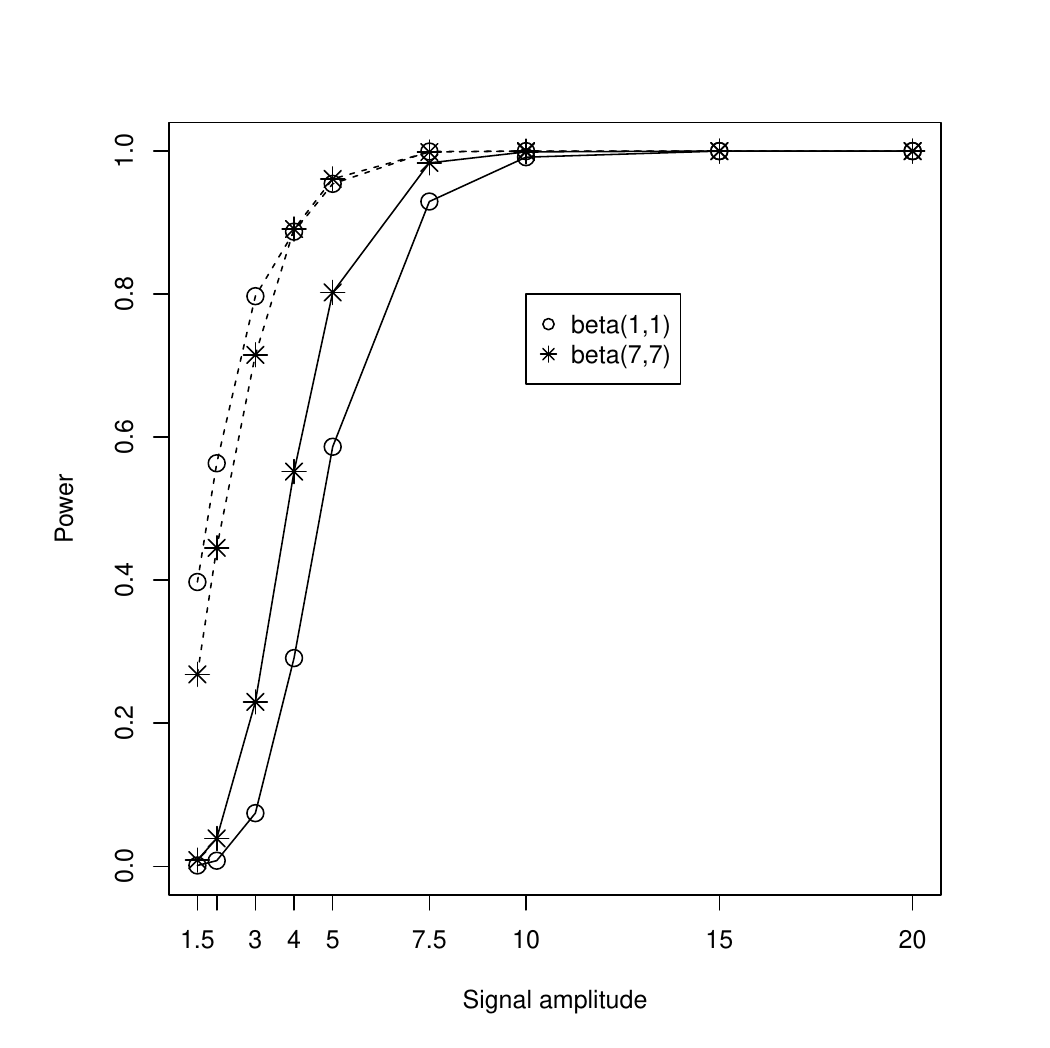}}
        \caption{Simulated data for 2-valued exchangeable covariates. Diffuse priors: Beta distributions with equal parameters}\label{fig:1}
   \end{figure}

\begin{figure}[!htbp]
        \centering
        \hbox{\hspace{-2cm}
        \includegraphics[width=0.7\textwidth]{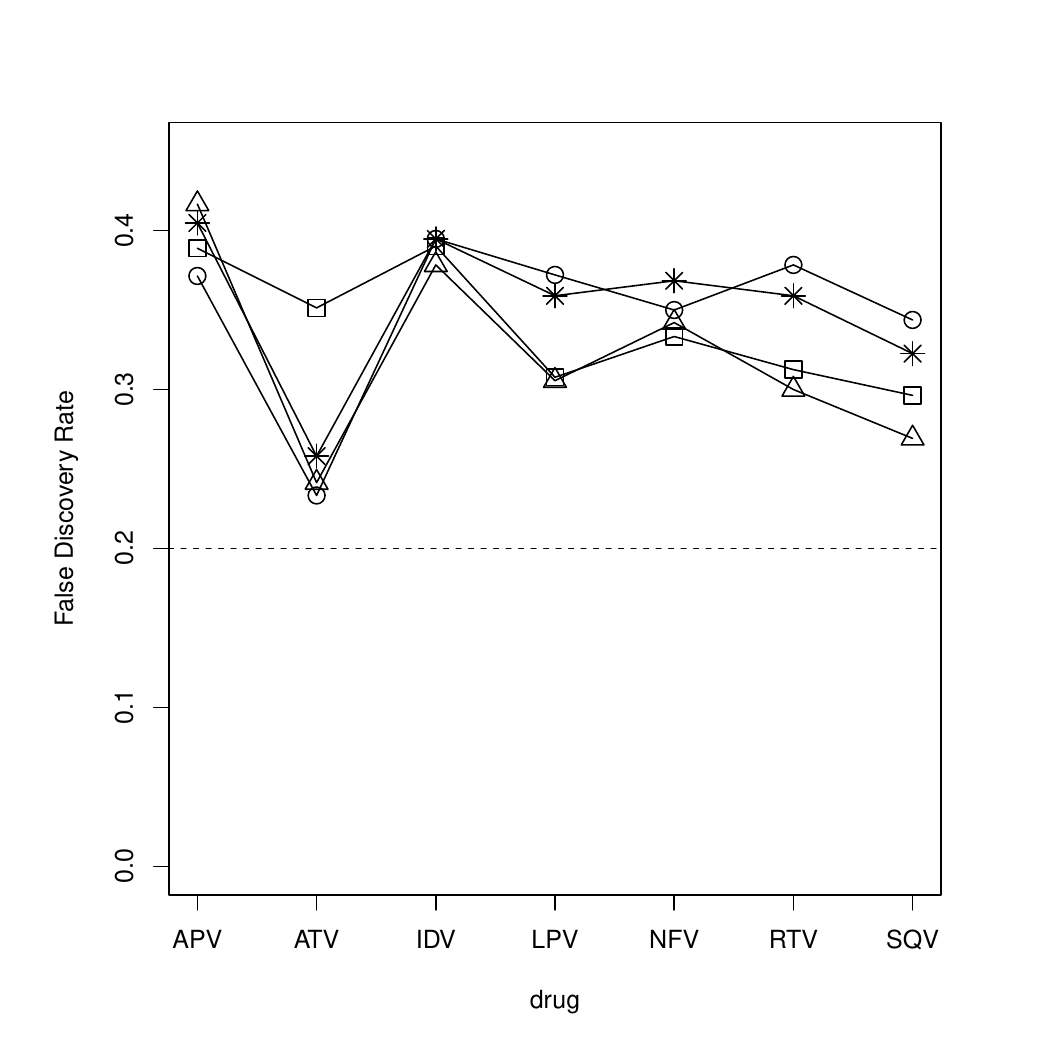}
        \hspace{-1cm}
        \includegraphics[width=0.7\textwidth]{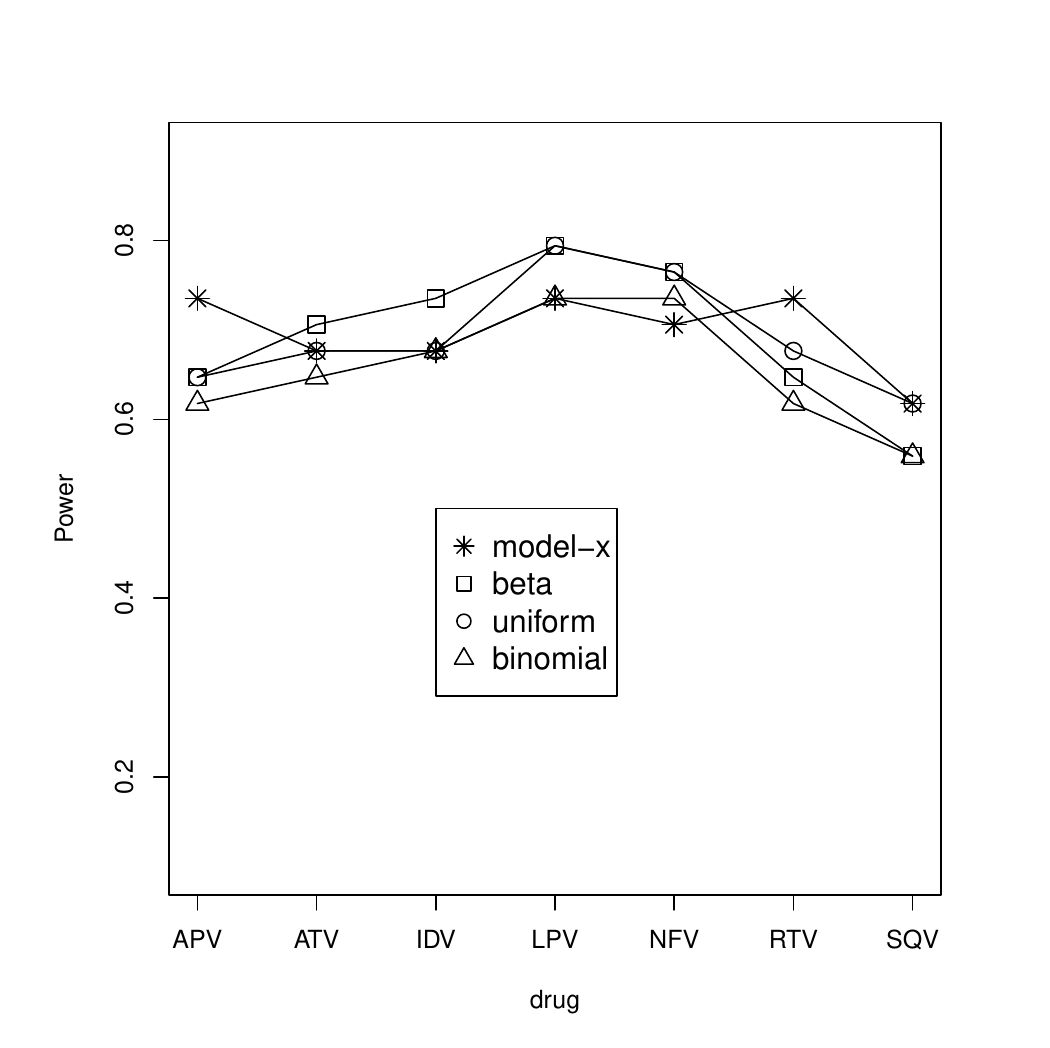}}
          \caption{Real data (HIV-1 dataset) for exchangeable 2-valued covariates. Priors: Beta, binomial and discrete uniform} \label{fig:3}
\end{figure}

\noindent In case of exchangeable binary covariates, the experiment has been performed for simulated and real data. Figures \ref{fig:1} and \ref{fig:2} correspond to simulated data while Figure \ref{fig:3} to real data (the HIV-1 dataset).

\medskip

\noindent In Figure \ref{fig:1}, the priors are beta distributions with equal parameters (i.e., $a=b$). It turns out that Model-$X$ is not able to keep the FDR under control. Indeed, the FDR of Model-$X$ systematically exceeds the nominal value of 0.1. On the contrary, the CIKs exhibit a very good behavior in terms of FDR. The situation is reversed as regards POW. The Model-$X$ is more powerful than the CIK, mainly for small values of the amplitude. However, if $a=b=7$ and the amplitude is $\ge 7.5$, the CIK and the Model-$X$ are essentially equivalent regarding the power. Moreover, as expected after Remark \ref{vg778nj}, the power of the CIK increases as the variance of the prior decreases.

\medskip

\begin{figure}[!htbp]
        \centering
        \hbox{\hspace{-2cm}\includegraphics[width=0.7\textwidth]{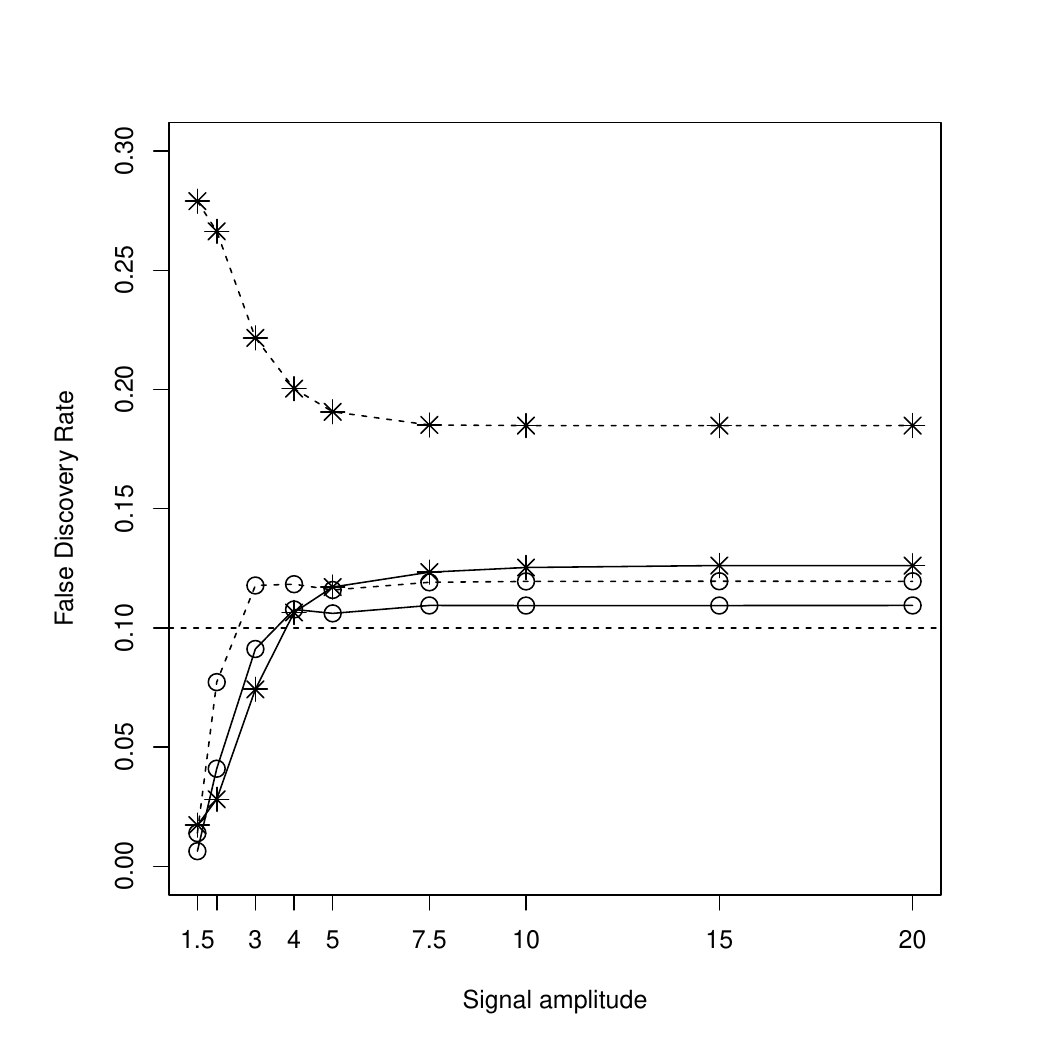}\hspace{-1cm}
         \includegraphics[width=0.7\textwidth]{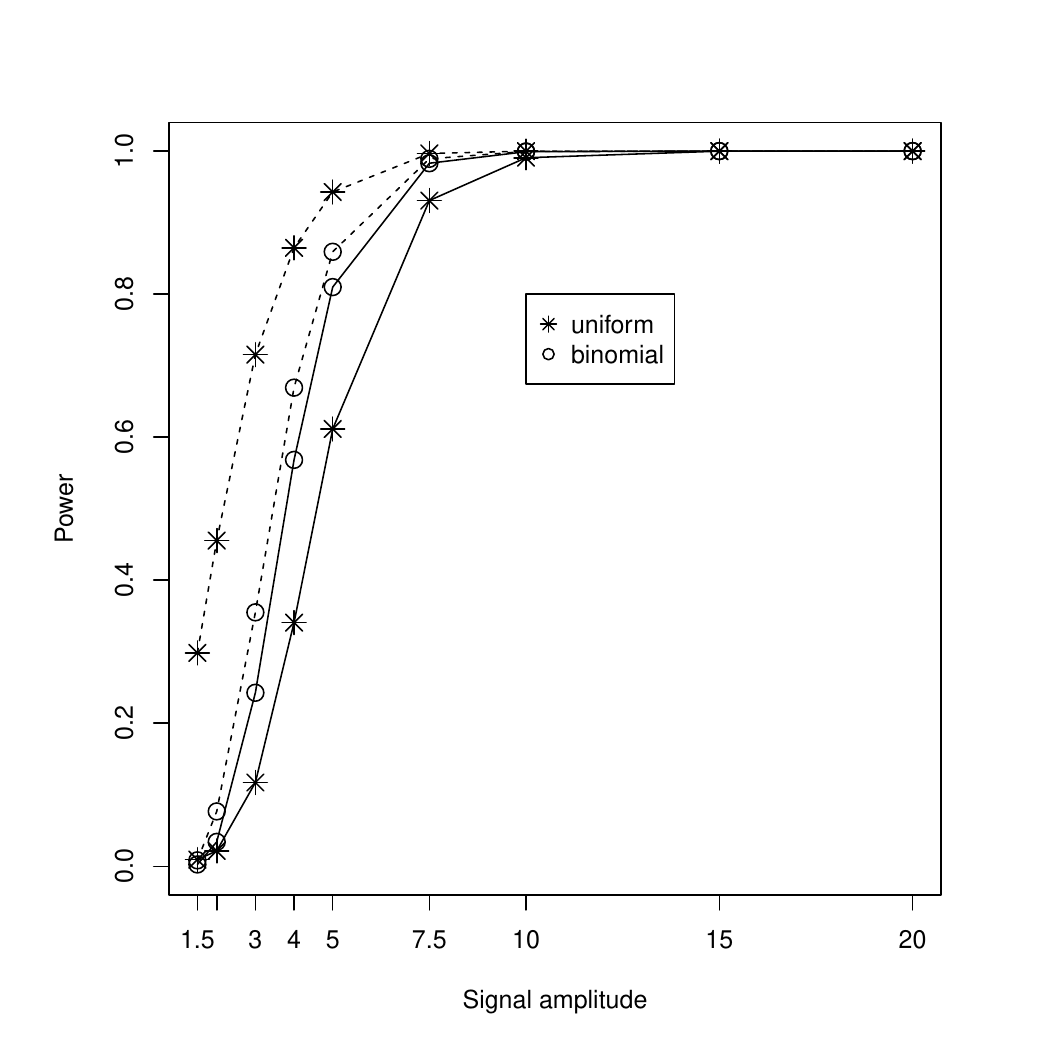}}
        \caption{Simulated data for exchangeable 2-valued covariates. Discrete priors: Uniform and binomial}\label{fig:2}
\label{fig:2}\end{figure}

\medskip

\noindent In Figure \ref{fig:2}, the priors are discrete: uniform and binomial. The situation is analogous to that of Figure \ref{fig:1}. Once again, the CIKs outperform Model-$X$ in terms of FDR but they are inferior in terms of POW. However, in case of binomial prior, the CIK and the Model-$X$ behave very similarly, as regards both FDR and POW.

\medskip

\medskip

\noindent Figure \ref{fig:3} refers to real data (the HIV-1 dataset). In this case, neither the CIKs nor the Model-$X$ keep the FDR under control. Apart from this, there are conflicting indications, without any clear pattern. Sometimes, the CIKs are better than Model-$X$ (with respect to both FDR and POW) but sometimes are worse.
\medskip

\begin{figure}[!htbp]
        \centering
        \hbox{\hspace{-2cm}\includegraphics[width=0.7\textwidth]{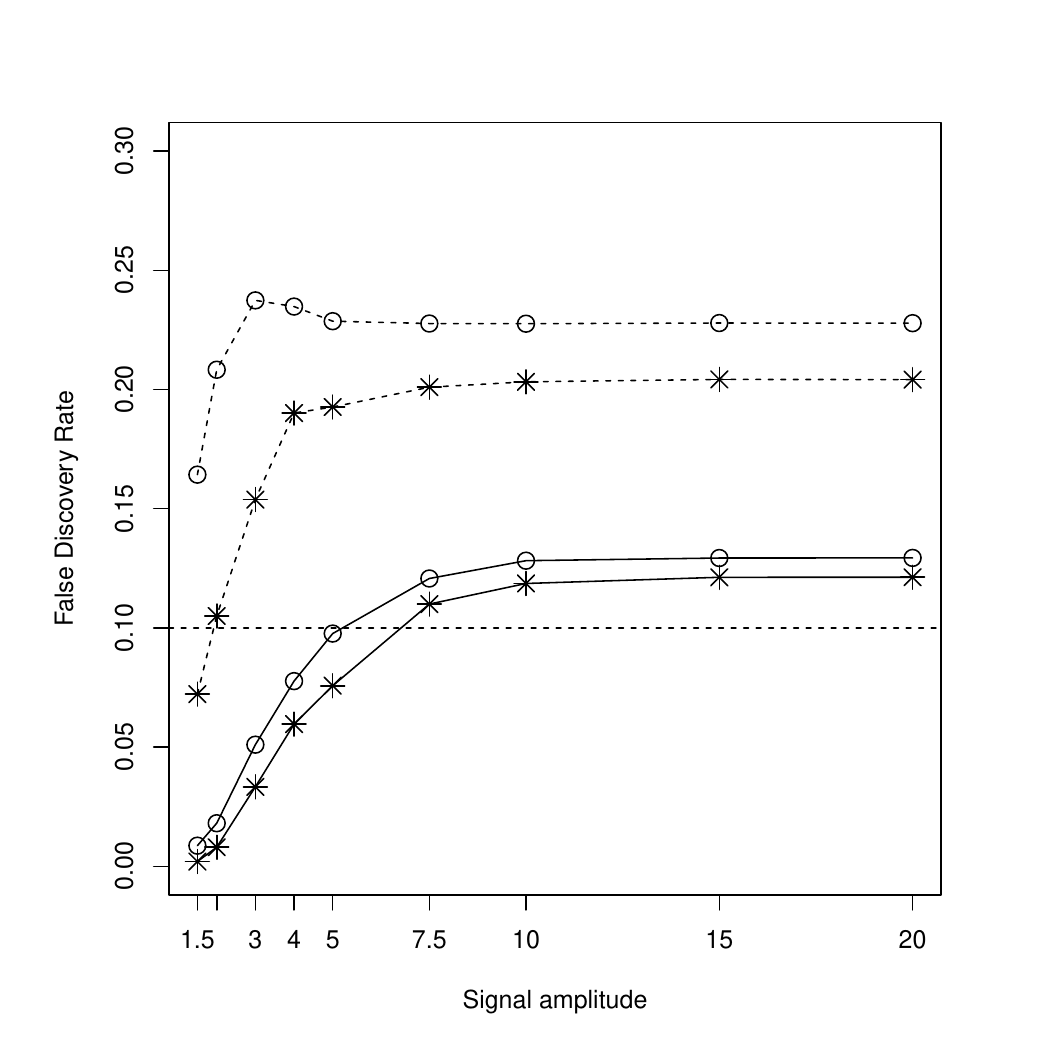}\hspace{-1cm}
         \includegraphics[width=0.7\textwidth]{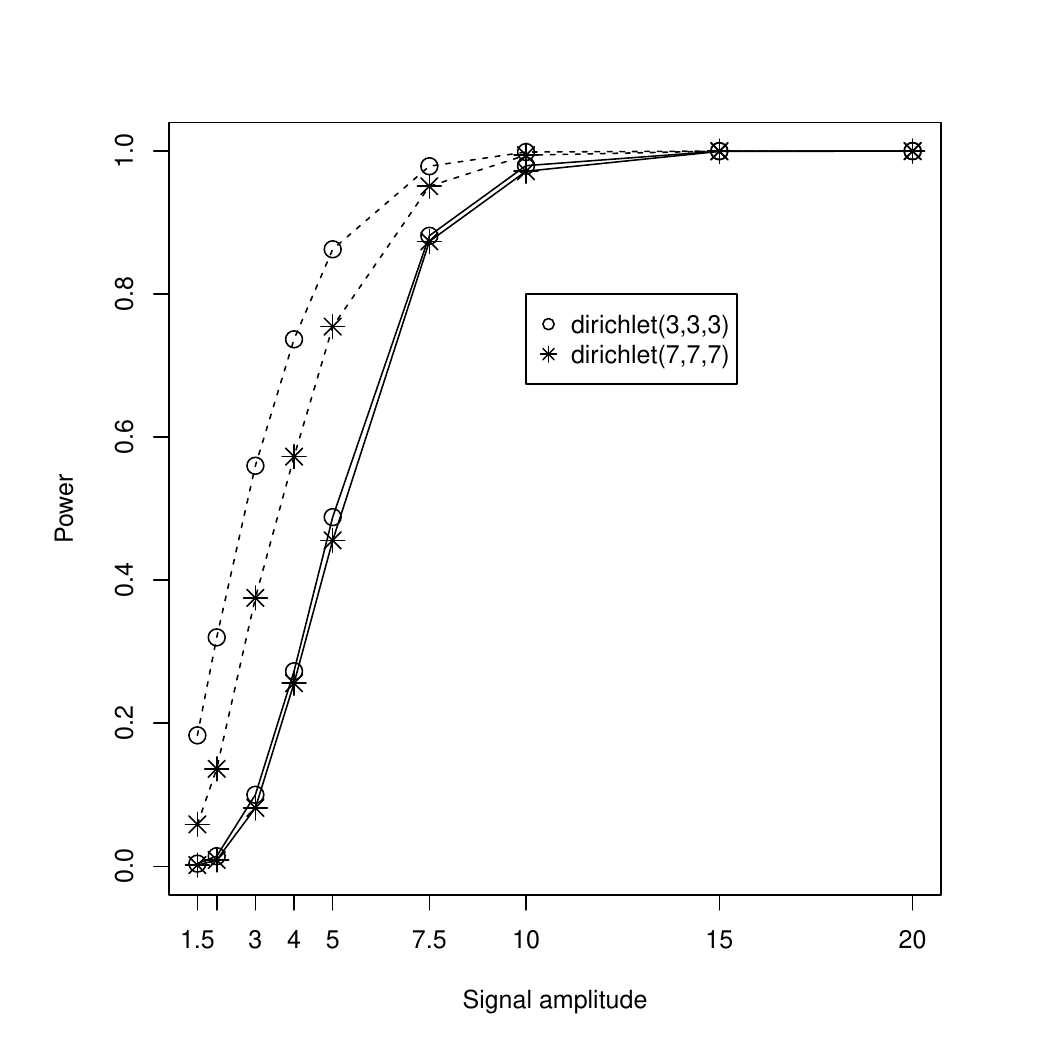}}
        \caption{Simulated data for exchangeable 3-valued covariates. Dirichlet priors with equal parameters}\label{fig:4}
   \end{figure}

\medskip

\noindent Finally, Figure \ref{fig:4} deals with simulated data for exchangeable 3-valued covariates. The priors are Dirichlet with equal parameters (i.e., $\alpha_0=\alpha_1=\alpha_2$ in the notation of Section \ref{b5rtg8n9}). Even in this case, the perfomances of both CIK and Model-$X$ get worse with respect to Figure \ref{fig:1}. That said, it seems still true that the CIKs are much better than Model-$X$ in terms of FDR and essentially equivalent in terms of POW for large amplitudes. However, they are inferior to Model-$X$ regarding POW for small amplitudes.

\medskip

\subsection{Partially exchangeable binary covariates} 
\begin{figure}[!htbp]
        \centering
        \hbox{
        \hspace{-2cm}\includegraphics[width=0.7\textwidth]{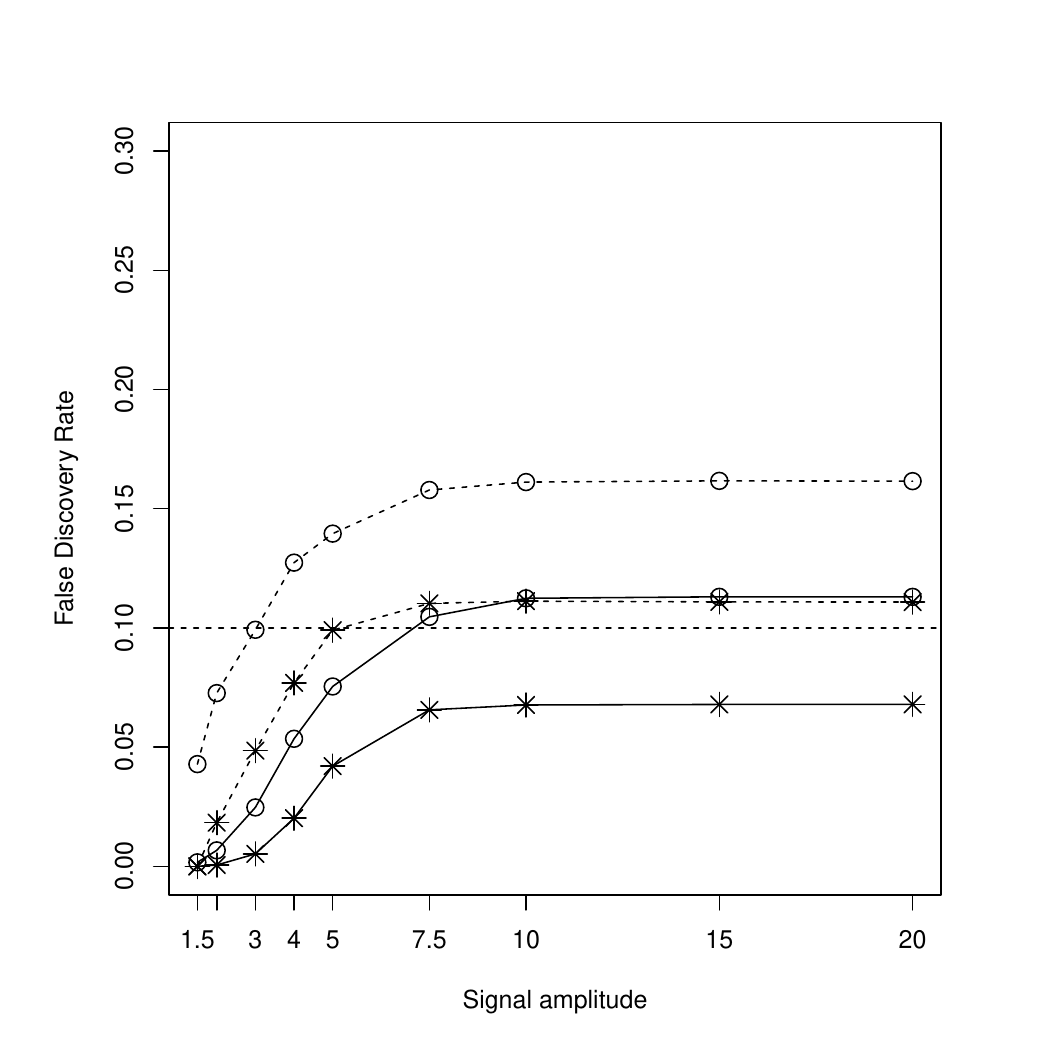}\hspace{-1cm}
         \includegraphics[width=0.7\textwidth]{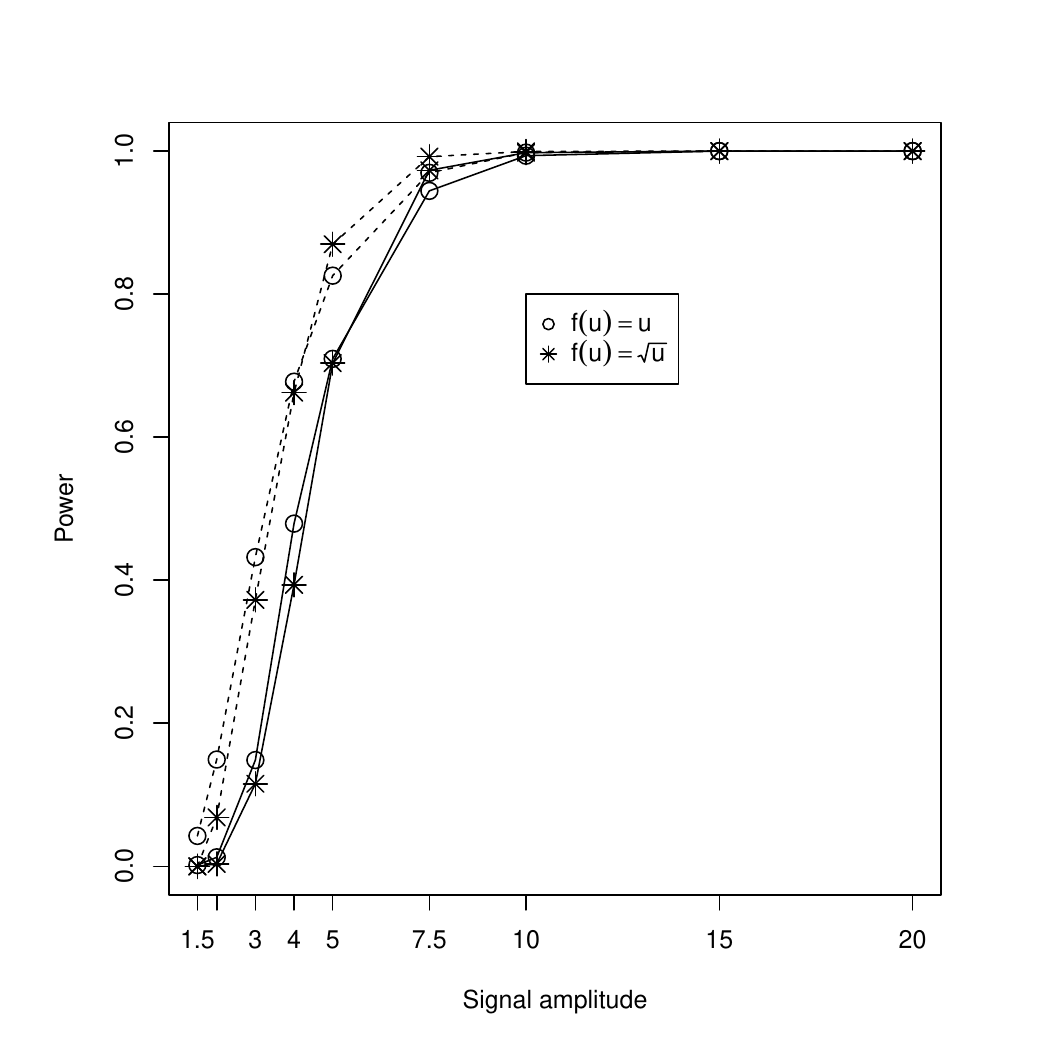}}
        \caption{Simulated data for partially exchangeable 2-valued covariates. Discrete priors with different choices of $f(u)$}\label{fig:5}
   \end{figure}

The simulation experiment for partially exchangeable binary covariates has been performed with $p=100$ and $k=50$. As shown in Figure \ref{fig:5}, the situation is still the same, namely, the CIKs are better than Model-$X$ regarding FDR but they are worse in terms of POW for small amplitudes. Incidentally, we note that the performances of the CIKs improve with respect to Figure \ref{fig:2} (exchangeable binary covariates with discrete priors).

\medskip

\medskip

\subsection{Concluding remarks} Drawing clear conclusions from the previous experiment is not trivial. While the CIKs are better than Model-$X$ regarding FDR, they are worse in terms of POW for small amplitudes. However, the good power of Model-$X$ could be probably expected; see \cite{WSBBC} and references therein. Hence, overall, the CIKs realize a reasonable ballance between FDR and POW; see also \cite[p. 139]{EJS}. In addition, unlike Model-$X$, the CIKs satisfy the basic constraint $P(X\in F^p)=1$.

\medskip

\noindent To close the paper, we briefly mention some hints for future research. They are essentially three, listed in increasing order of importance. The first is looking for further priors. This is not very urgent, however. Even if quite standard, the priors suggested in this paper are flexible enough to cover most real situations. Moreover, the main purpose should be to make KP effective, not to introduce new sophisticated priors. A second hint is to drop the exchangeability assumption. As argued in the Introduction, in case of categorical covariates, exchangeability is often a reasonable condition. But of course this in not always true. In several problems, partial exchangeability (possibly with more than two groups) is more suitable. We would not go beyond partial exchangeability, since otherwise the theory of Section \ref{vh87} does not apply and building knockoffs becomes very hard. Finally, and more importantly, the ``second step of KP" should be revised and adapted to categorical covariates. By the second step of KP, we mean that step of the procedure where the knockoff $\widetilde{X}$ has been already built, and one has to compare it with $X$ and $Y$ in order to decide which covariates are to be discarded. In this paper, this important step has been neglected. In particular, in the numerical experiment, we exploited the R-cran package \verb"knockoff" which is built having quantitative covariates in mind. Most probably, in case of categorical covariates, this second step should be rethought. Something similar has been done in \cite{STATMED}.

\end{document}